\documentclass[1p]{elsarticle}
\usepackage[utf8]{inputenc}
\usepackage[UKenglish]{babel}
\usepackage{amsmath}
\usepackage{amsthm}
\usepackage{amscd}
\usepackage{amssymb}
\usepackage{MnSymbol}
\usepackage{latexsym}
\usepackage{eucal}
\usepackage{dsfont}
\usepackage{mathtools}
\usepackage{enumitem}

\usepackage[colorlinks,pdftex]{hyperref}
\hypersetup{
linkcolor=black,
citecolor=black,
pdftitle={On the domain of fractional Laplacians and related generators of Feller processes},
pdfauthor={Franziska K\"uhn, Ren\'e L. Schilling},
pdfkeywords={},
}

\newcommand\specialsection[1]{%
   \let\savesection\thesection
   \renewcommand*{\thesection}{\Alph{section}}}
   
\newcommand\defaultsection{%
   \renewcommand\thesection{\savesection}}

\newtheorem{theorem}{Theorem}[section]
\newtheorem{lemma}[theorem]{Lemma}
\newtheorem{corollary}[theorem]{Corollary}
\newtheorem{proposition}[theorem]{Proposition}

\newdefinition{remark}[theorem]{Remark}
\newdefinition{example}[theorem]{Example}
\newdefinition{definition}[theorem]{Definition}

\DeclareMathOperator \re {Re}
\DeclareMathOperator \im {Im}

\DeclareMathOperator \spt {supp}

\DeclareMathOperator \lip {Lip}

\DeclareMathOperator \tr {tr}

\newcommand{\I}{\mathds{1}}

\newcommand\fa{\qquad \text{for all \ }}

\newcommand{\cadlag}{c\`adl\`ag }

\newcommand\mc[1] {\mathcal{#1}}
\newcommand\mbb[1] {\mathds{#1}}

\newcommand{\eps}{\varepsilon}

\newcommand{\ball}[2]{B(#1,#2)} 
\newcommand{\coball}[2]{B^c(#1,#2)} 
\newcommand{\cball}[2]{B[#1,#2]} 
\newcommand{\rd}{{\mbb{R}^d}}

\makeatletter
\def\mov@rlay#1#2{\leavevmode\vtop{%
   \baselineskip\z@skip \lineskiplimit-\maxdimen
   \ialign{\hfil$\m@th#1##$\hfil\cr#2\crcr}}}
\newcommand{\charfusion}[3][\mathord]{
    #1{\ifx#1\mathop\vphantom{#2}\fi
        \mathpalette\mov@rlay{#2\cr#3}
      }
    \ifx#1\mathop\expandafter\displaylimits\fi}

\newcommand{\intx}{\charfusion[\mathop]{\int}{\innerx}}

\newcommand{\innerx}{%
  \mathchoice
    {\mbox{\fontsize\f@size{0}\normalfont\hspace{-0.15pt}$\times$}}
    {\mbox{\fontsize\sf@size{0}\normalfont\hspace{-0.9pt}$\times$}}
    {\mbox{\fontsize\ssf@size{0}\normalfont$\times$}}
    {\mbox{\tiny\normalfont\ $\times$}}%
}
\makeatother

\makeatletter
\def\@fnsymbol#1{\ensuremath{\ifcase#1\or *\or **\or \dagger\or \ddagger\or
   \mathsection\or \mathparagraph\or \|\or \dagger\dagger
   \or \ddagger\ddagger \else\@ctrerr\fi}}

\newcommand{\lint}{\intx} 

\begin{document}

\title{On the domain of fractional Laplacians and related generators of Feller processes}	
	
\author[FK]{Franziska K\"{u}hn\corref{cor1}}
\cortext[cor1]{Corresponding author}
\ead{franziska.kuhn@math.univ-toulouse.fr}
\address[FK]{Institut de Math\'ematiques de Toulouse, Universit\'e Paul Sabatier III Toulouse, 118 Route de Narbonne, 31062 Toulouse, France. \emph{On leave from:} TU Dresden, Fachrichtung Mathematik, Institut f\"{u}r Mathematische Stochastik, 01062 Dresden, Germany}

\author[RS]{Ren\'e L.\ Schilling}
\ead{rene.schilling@tu-dresden.de}
\address[RS]{TU Dresden, Fakult\"{a}t Mathematik, Institut f\"{u}r Mathematische Stochastik, 01062 Dresden, Germany}


\begin{abstract}
    In this paper we study the domain of the generator of stable processes, stable-like processes and more general pseudo- and integro-differential operators which naturally arise both in analysis and as infinitesimal generators of L\'evy- and L\'evy-type (Feller) processes. In particular we obtain conditions on the symbol of the operator ensuring that certain (variable order) H\"{o}lder and H\"{o}lder-Zygmund spaces are in the domain. We use tools from probability theory to investigate the small-time asymptotics of the generalized moments of a L\'evy or L\'evy-type process $(X_t)_{t \geq 0}$,
    \begin{equation*}
    	\lim_{t \to 0} \frac 1t\left(\mathbb{E}^x f(X_t)-f(x)\right), \quad x\in\mathbb{R}^d,
    \end{equation*}
    for functions $f$ which are not necessarily bounded or differentiable. The pointwise limit exists for fixed $x \in \mathbb{R}^d$ if $f$ satisfies a H\"{o}lder condition at $x$. Moreover, we give sufficient conditions which ensure that the limit exists uniformly in the space of continuous functions vanishing at infinity. As an application we prove that the domain of the generator of $(X_t)_{t \geq 0}$ contains certain H\"{o}lder spaces of variable order. Our results apply, in particular, to stable-like processes, relativistic stable-like processes, solutions of L\'evy-driven SDEs and L\'evy processes.
 \end{abstract}

\maketitle

\section{Introduction} \label{intro}

Since the pioneering work of Caffarelli and Silvestre on fractional powers of the Laplacian, see \cite{silvestre,caf-sil}, a lot of work has been devoted to fractional powers of the Laplacian from the analytical point of view, we refer to \cite{CSI,CSII,caf-sil14,grubb,ros-oton} to mention but a few.

The fractional power of the Laplacian is also the generator of a stochastic process with stationary and independent increments (a L\'evy process), which allows us to use probabilistic methods for its investigation. In fact, fractional powers of the Laplacian are just a special case of generators of L\'evy processes and -- if one allows for generators with variable coefficients -- of the more general class of Feller processes, the classic result is \cite{courrege}, see \cite{ltp} for a recent survey. Over the past two and a half decades these operators have been studied from both the analytical community but most of all the probability community, see \cite{bogdan,chen-kim-song,chen-kum,jacob123,kas-mim,sato,taira}.

Of particular importance is a good understanding of the domain of these operators which, in general, have a representation as pseudo-differential as well as integro-differential operator. This is partly due to the fact that for elements in their domains we can construct interesting martingales.

In this paper we study in great detail the domains of rather general generators of Feller processes and, by using probabilistic techniques in combination with analytic techniques, we succeed in finding precise conditions in terms of (variable-order) H\"{o}lder and Lipschitz function spaces to belong to these domains, see Theorem \ref{app-3} (for L\'evy processes and generators with constant coefficients) and Theorem \ref{app-5} (for Feller processes and generators with variable coefficients). As far as we are aware, these results extend known results for fractional powers of the Laplacian (including those of variable order of differentiability).

For a $d$-dimensional L\'evy process $(L_t)_{t \geq 0}$ with L\'evy triplet $(b,Q,\nu)$ the family of measures $(p_t)_{t>0}$ on $(\mbb{R}^d \setminus \{0\},\mc{B}(\mbb{R}^d \setminus \{0\})$ defined by
\begin{equation*}
	p_t(B) := \frac 1t \mbb{P}(L_t \in B), \qquad t>0,\; B \in \mc{B}(\mbb{R}^d \setminus \{0\})
\end{equation*}
converges vaguely to the L\'evy measure $\nu$, i.\,e.\
\begin{equation}\label{levy-pw}
	\lim_{t \to 0} \frac{1}{t} \mbb{E}f(L_t) = \int_{\rd\setminus\{0\}}  f(y) \, \nu(dy)
\end{equation}
holds for any continuous function $f$ with compact support in $\mbb{R}^d \setminus \{0\}$, cf.\ \cite[Lemma 2.16]{ltp} or \cite[Proposition 18.2]{berg}. By the portmanteau theorem, this implies the following small-time asymptotics
\begin{equation}\label{levy-ihke}
	\lim_{t \to 0} \frac{1}{t} \mbb{P}(L_t \in B) = \nu(B)
\end{equation}
for any Borel set $B \in \mc{B}(\mbb{R}^d \setminus \{0\})$ such that $0\notin\bar B$ and the topological boundary $\partial B$ is a $\nu$-null set.  Jacod \cite{jac} proved that the small-time asymptotics \eqref{levy-pw} extends to continuous bounded functions $f: \mbb{R}^d \to \mbb{R}$ with $f(0)=0$ which satisfy a H\"older condition at $x=0$,
\begin{equation*}
	|f(x)-f(0)| \leq |x|^{\alpha} \fa |x| \leq 1
\end{equation*}
where $\alpha \in (0,2)$ is a suitable constant depending on the L\'evy triplet $(b,Q,\nu)$, see \cite{jac} or \cite[p.~2]{fig} for details. More recently, Figueroa--L\'opez \cite{fig} showed that the assumption on the boundedness of $f$ can be replaced by a much weaker integrability condition which basically ensures that the expectation $\mbb{E}f(L_t)$ is exists for any $t>0$.

In the first part of this paper, Section~\ref{s-erg}, we establish similar results for the class of L\'evy-type processes which includes, in particular, L\'evy processes, affine processes, solutions of L\'evy-driven stochastic differential equations, and stable-like processes. We will show that any L\'evy-type process $(X_t)_{t \geq 0}$ with rich domain and characteristics $(b(x),Q(x),\nu(x,dy))$ satisfies \begin{equation}\label{ltp-ihke}
	\lim_{t \to 0} \frac{1}{t} \mbb{P}^x(X_t-x \in B) = \nu(x,B) \fa  x \in \mbb{R}^d
\end{equation}
which is the analogue of \eqref{levy-ihke}, cf.\ Corollary~\ref{erg-7}; again $B \in \mc{B}(\mbb{R}^d \setminus \{0\})$ is a Borel set such that $0\notin\bar B$ and $\nu(x,\partial B)=0$. Because of the small-time asymptotics \eqref{ltp-ihke}, we have for fixed $x \in \mbb{R}^d$ \begin{equation*}
	\lim_{t \to 0} \frac{1}{t} \left(\mbb{E}^x f(X_t)-f(x)\right) = \int_{\rd\setminus\{0\}}  (f(x+y)-f(x)) \, \nu(x,dy)
\end{equation*}
for any continuous function $f$ with compact support in $\mbb{R}^d \setminus \{x\}$. Using a localized version of a maximal inequality, cf.\ Lemma~\ref{erg-3}, we will show that for a rich L\'evy-type process $(X_t)_{t \geq 0}$ and fixed $x \in \mbb{R}^d$ the pointwise limit
\begin{equation}\label{ltp-pw}
	\lim_{t \to 0} \frac{1}{t} \left(\mbb{E}^x f(X_t)-f(x)\right)
\end{equation}
exists for a much larger class of functions. More precisely, we will establish the small-time asymptotics \eqref{ltp-pw} for functions $f:\mbb{R}^d \to \mbb{R}$ which satisfy a H\"older condition at $x$, cf.\ Theorem~\ref{erg-9} and~\ref{erg-10}, and need not be bounded, see Theorem~\ref{erg-13}.

In the second part, Section~\ref{s-app}, we turn to the question under which assumptions on a continuous function $f$ vanishing at infinity -- we write $f \in C_0(\mbb{R}^d)$ for short -- the limit
\begin{equation}\label{ltp-uni}
	\lim_{t \to 0} \frac{1}{t} \left(\mbb{E}^x f(X_t)-f(x)\right)
\end{equation}
exists uniformly (in $x$) for a rich L\'evy-type process $(X_t)_{t \geq 0}$ with bounded coefficients. This is equivalent to asking for sufficient conditions which ensure that a function $f \in C_0(\mbb{R}^d)$ is contained in the domain $\mc{D}(A)$ of the generator $A$ of $(X_t)_{t \geq 0}$. The main results in Section~\ref{s-app} are Corollary~\ref{app-11} and Corollary~\ref{app-13} which state that $\mc{D}(A)$ contains certain H\"older spaces of variable order. Our results apply, in particular, to L\'evy processes, cf.\ Theorem~\ref{app-3}; for instance, if $(L_t)_{t \geq 0}$ is an isotropic $\alpha$-stable L\'evy process, $\alpha \in (0,1)$, then the H\"older space $\mc{C}_0^{\beta}$ -- see \eqref{def-eq2} below for a precise definition -- is contained in the domain of the generator $A$ of $(L_t)_{t \geq 0}$ for any $\beta \in (\alpha,1]$, and we have
\begin{equation*}
	Af(x) = \int_{\rd\setminus\{0\}}  (f(x+y)-f(x)) \, \nu(dy), \qquad f \in \mc{C}_0^{\beta},\; x \in \mbb{R}^d.
\end{equation*}
At the end of Section~\ref{s-app} we discuss several examples, including stable-like dominated processes (Example~\ref{app-15}), solutions of L\'evy-driven SDEs (Example~\ref{app-19}), stable-like processes and relativistic stable-like processes (Example~\ref{app-17}).

\section{Basic definitions and notation} \label{def}

We consider the Euclidean space $\mbb{R}^d$ with the canonical scalar product $x \cdot y := \sum_{j=1}^d x_j y_j$ and the Borel $\sigma$-algebra $\mc{B}(\mbb{R}^d)$ generated by the open balls $\ball{x}{r} := \{y \in \mbb{R}^d; |y-x|<r\}$ and closed balls $\cball{x}{r} := \{y \in \mbb{R}^d; |y-x| \leq r\}$. We write $\spt f$ for the support of a function $f: \mbb{R}^{n} \to \mbb{R}^d$ and $\{f \in B\} = f^{-1}(B)$ denotes the preimage of a set $B \subseteq \mbb{R}^d$ under $f$. A function $f: \mbb{R}^d \to [0,\infty)$ is called \emph{submultiplicative} if there exists a constant $c>0$ such that
\begin{equation}\label{def-eq0}
	f(x+y) \leq c f(x) f(y) \fa x,y \in \mbb{R}^d.
\end{equation}
Later on we will use that submultiplicative functions grow at most exponentially, cf.\ \cite[Lemma 25.5]{sato}. For a set $B \subseteq \mbb{R}^d$ we use $\partial B$ to denote the topological boundary of $B$. We use $\lint$ and $\lint_B$ as a shorthand for $\int_{\mbb{R}^d \setminus \{0\}} $ and $\int_{B\setminus\{0\}}$, respectively.

\medskip
\textbf{Function spaces:} The smooth functions with compact support are denoted by $C_c^{\infty}(\mbb{R}^d)$, and $C_0(\mbb{R}^d)$ is the space of continuous functions $f: \mbb{R}^d \to \mbb{R}$ vanishing at infinity. Superscripts $k\in\mbb{N}$ are used to denote the order of differentiability, e.\,g.\ $f \in C_0^k(\mbb{R}^d)$ means that $f$ and its derivatives up to and including order $k$ are $C_0(\mbb{R}^d)$-functions. We define \emph{H\"{o}lder spaces} by
	\begin{align} \label{def-eq2} \begin{aligned}
	\mc{C}^{\alpha}
	&:= \left\{f \in C_0(\mbb{R}^d); \: \|f\|_{\alpha} := \sup_{x,y \in \mbb{R}^d} \frac{|f(x)-f(y)|}{|x-y|^{\alpha}}<\infty \right\}, \quad \alpha \in [0,1] \\
	\mc{C}_0^{1,\alpha}
	&:= \left\{f \in C_0^1(\mbb{R}^d); \: \nabla f \in \mc{C}_0^{\alpha} \right\}, \quad \alpha \in [0,1] \end{aligned}
\end{align}
Since there are various concepts of H\"{o}lder (or Lipschitz) spaces in the literature, let us explain the relations to the other function spaces. There are the ``classical'' H\"{o}lder spaces $C^{\alpha}$ equipped with the norm
	\begin{equation}
	\sum_{j=0}^{\lfloor \alpha \rfloor}
	\sum_{\substack{\beta \in \mbb{N}_0^d \\ |\beta| = j}} \|\partial^{\beta} f\|_{\infty}
	+   \max_{\substack{\beta \in \mbb{N}_0^d \\ |\beta| = \lfloor \alpha \rfloor}}   \sup_{x \neq y} \frac{|\partial^{\beta} f(x)-\partial^{\beta} f(y)|}{|x-y|^{\alpha-\lfloor \alpha \rfloor}} \label{app-st1} \tag{$\star$}
	\end{equation}
	where $\lfloor \alpha \rfloor$ denotes the biggest natural number less or equal than $\alpha$. On the other hand, there are the Zygmund--H\"{o}lder spaces $\mathfrak{C}^{\alpha}$ consisting of all functions $f \in C^k$ such that the norm
	\begin{equation*}
		\sum_{j=0}^k
		\sum_{\substack{\beta \in \mbb{N}_0^d \\ |\beta| = j}}
		\|\partial^{\beta} f\|_{\infty}
		+
		\max_{\substack{\beta \in \mbb{N}_0^d \\ |\beta| = k}}
		\sup_{\substack{x,h \in \mbb{R}^d \\ h \neq 0}}
		\frac{|\partial^{\beta} f(x+h)+\partial^{\beta} f(x-h)-2\partial^{\beta}f(x)|}{|h|^{s}}
	\end{equation*}
	is finite where $s \in (0,1]$ and $k \in \mbb{N}$ are chosen such that $\alpha = k+s$, see Triebel \cite[pp.~34]{triebel}. If $\alpha \in (0,\infty)\setminus \mbb{N}$ then $\mathfrak{C}^{\alpha} = C^{\alpha}$, cf.\ \cite[Theorem 1(b), p.~201]{triebel-inter}; however for $\alpha \in \mbb{N}$ we have a strict inclusion: $\mathfrak{C}^{\alpha} \supsetneq C^{\alpha}$. For $\alpha=1$ it is possible to show that $\mathfrak{C}^1$ is strictly larger than the space of Lipschitz continuous functions $\lip$ (cf.\ \cite[p.~148]{stein}) which is, in turn, strictly larger than $C^1$.
	 There are the following relations between the H\"{o}lder spaces introduced in \eqref{def-eq2} and the just mentioned function spaces:
	\begin{align*}
		\mc{C}_0^{\alpha} &= C^{\alpha} \cap C_0(\mbb{R}^d) = \mathfrak{C}^{\alpha} \cap C_0(\mbb{R}^d),\qquad\alpha \in (0,1), \\
		\mc{C}_0^{1,\alpha-1} &= C^{\alpha} \cap C_0^1(\mbb{R}^d) = \mathfrak{C}^{\alpha} \cap C_0^1(\mbb{R}^d),\qquad\alpha \in (1,2)
	\end{align*}
	and
	\begin{equation*}
		\mc{C}_0^1 = \lip \cap C_0(\mbb{R}^d), \qquad
		\mc{C}_0^{1,0} = C^1 \cap C_0(\mbb{R}^d).
	\end{equation*}

\textbf{L\'evy(-type) Processes:} Throughout, $(\Omega,\mc{A},\mbb{P})$ denotes a probability space. A stochastic process $(L_t)_{t \geq 0}$ is called a \emph{L\'evy process} if it has stationary and independent increments, $L_0=0$ almost surely and the sample paths $t \mapsto L_t(\omega)$ are \cadlag (right-continuous with finite left-hand limits) for almost all $\omega \in \Omega$. By the L\'evy-Khintchine formula, every L\'evy process can be uniquely characterized by its \emph{characteristic exponent} $\psi(\xi) := -\log\mbb{E} e^{i\xi\cdot X_1}$,
\begin{equation}\label{def-eq1}
    \psi(\xi)
    = -ib \cdot \xi + \frac{1}{2} \xi \cdot Q \xi
    + \lint  \left(1-e^{iy \cdot \xi}+iy \cdot \xi \I_{(0,1)}(|y|)\right)  \nu(dy),
    \quad \xi \in \mbb{R}^d,
\end{equation}
where $(b,Q,\nu)$ is the \emph{L\'evy triplet} consisting of the drift $b \in \mbb{R}^d$, the symmetric positive semidefinite diffusion matrix $Q \in \mbb{R}^{d \times d}$ and the L\'evy measure $\nu$ on $(\mbb{R}^d \setminus \{0\},\mc{B}(\mbb{R}^d \setminus \{0\}))$ satisfying $\lint \min\{|y|^2,1\} \, \nu(dy)<\infty$. A function $\psi: \mbb{R}^d \to \mbb{C}$ with $\psi(0)=0$ is called \emph{continuous negative definite} if it admits a L\'evy--Khintchine representation of the form \eqref{def-eq1}.

A \emph{L\'evy-type process} is a Markov process whose transition semigroup is a Feller semigroup; for further details see e.\,g.\ \cite{ltp}. Without loss of generality, we may assume that the sample paths of a L\'evy-type process are c\`adl\`ag. If $C_c^{\infty}(\mbb{R}^d)$ is contained in the domain $\mc{D}(A)$ of the generator $A$ of a L\'evy-type process $(X_t)_{t \geq 0}$, then we call $(X_t)_{t \geq 0}$ a \emph{rich} L\'evy-type process. L\'evy-type processes are also known as Feller processes, and we will use both terms synonymously. Our main reference for Feller processes is the monograph \cite{ltp}. If $(X_t)_{t \geq 0}$ is a rich L\'evy-type process with generator $A$, then $A|_{C_c^{\infty}(\mbb{R}^d)}$ is a \emph{pseudo-differential operator},
\begin{equation*}
    Af(x)
    =-q(x,D)f(x)
    := - \int_{\mbb{R}^d} e^{i \, x \cdot \xi} q(x,\xi) \widehat{f}(\xi) \, d\xi, \qquad f \in C_c^{\infty}(\mbb{R}^d),\; x \in \mbb{R}^d
\end{equation*}
where $\widehat{f}(\xi) := (2\pi)^{-d} \int_{\mbb{R}^d} e^{-ix \cdot \xi} f(x) \, dx$ denotes the Fourier transform of $f$ and
\begin{equation}\label{def-eq3}
	q(x,\xi)
    = q(x,0) - i b(x) \cdot \xi + \frac{1}{2} \xi \cdot Q(x) \xi
    + \lint  \left(1-e^{i y \cdot \xi}+ i y\cdot \xi \I_{(0,1)}(|y|)\right)  \nu(x,dy)
\end{equation}
is the negative definite \emph{symbol}, cf.\ \cite[Theorem 2.21]{ltp}.  By \cite[Theorem 2.30]{ltp}, continuity of $x\mapsto q(x,0)$ implies that the mapping $x \mapsto q(x,\xi)$ is continuous for all $\xi \in \mbb{R}^d$. Probabilistically, the term $q(x,0)$ leads to a(n exponential) killing of the process, while analytically it acts like a multiplication operator. Both cases are not interesting for our study and we will assume from now on that $q(x,0)=0$. 

For each fixed $x \in \mbb{R}^d$ the tuple $(b(x),Q(x),\nu(x,dy))$ is a L\'evy triplet. We call the family $(b(x),Q(x),\nu(x,dy))_{x \in \mbb{R}^d}$ the \emph{characteristics} of $q$ and use $(b,Q,\nu)$ as a shorthand. It is not difficult to see that
\begin{align*}
    Af(x)
    &= b(x) \cdot \nabla f(x)+\frac{1}{2} \tr\left(Q(x) \cdot \nabla^2 f(x)\right) \\
    &\quad + \lint  \left(f(x+y)-f(x)-\nabla f(x) \cdot y \I_{(0,1)}(|y|)\right)  \nu(x,dy)
\end{align*}
for any $f \in C_c^{\infty}(\mbb{R}^d)$, see e.\,g.\ \cite[Theorem 2.21]{ltp}, where $\nabla^2 f$ denotes the Hessian and $\tr A$ the trace of a matrix $A$.   We say that a rich L\'evy-type process $(X_t)_{t \geq 0}$ has \emph{bounded coefficients} if its symbol $q$ has bounded coefficients, i.\,e.\ there exists a constant $c>0$ such that $|q(x,\xi)| \leq c(1+|\xi|^2)$ for all $x,\xi \in \mbb{R}^d$. We will frequently use the following result from \cite[Proposition 2.27(d), Theorem 2.31]{ltp}.

\begin{theorem} \label{def-3}
	Let $q$ be given by \eqref{def-eq3} such that $q(x,0)=0$. For any compact set $K \subseteq \mbb{R}^d$:
    \begin{enumerate}
    	\item $C_K := \sup_{x \in K} \sup_{|\xi| \leq 1} |q(x,\xi)| <\infty$,
    	\item $\sup_{x \in K} |q(x,\xi)| \leq 2 C_K (1+|\xi|^2)$ for all $\xi \in \mbb{R}^d$,
    	\item $\sup_{x \in K} (|b(x)|+|Q(x)|+ \lint  (|y|^2 \wedge 1) \, \nu(x,dy)) < \infty$.
    \end{enumerate}
    If $q$ has bounded coefficients, then the statements also hold for $K=\mbb{R}^d$.
\end{theorem}

We define, following \cite{rs97}, for fixed $x_0 \in \mbb{R}^d$ the generalized Blumenthal--Getoor index at $\infty$
\begin{equation}
    \beta_{\infty}^{x_0} := \inf \left\{\gamma>0; \lim_{r \to \infty} \frac{1}{r^{\gamma}} \sup_{|\xi| \leq r} |q(x_0,\xi)|< \infty \right\}. \label{erg-eq4}
\end{equation}
Since any continuous negative definite function grows at most quadratically at infinity, we have $\beta_{\infty}^{x_0} \in [0,2]$ for any $x_0 \in \mbb{R}^d$; moreover,
\begin{equation}\label{erg-eq6}
\lint_{|y| \leq 1} |y|^{\beta} \, \nu(x_0,dy)< \infty \fa \beta > \beta_{x_0}^{\infty}.
\end{equation}
If $q(x_0,\cdot)$ has no diffusion part, i.\,e.\ $Q(x_0)=0$, and satisfies the sector condition, i.\,e.\ if there exists a constant $C>0$ such that $|\im q(x_0,\xi)| \leq C \re q(x_0,\xi)$ for all $\xi \in \mbb{R}^d$, then
\begin{equation}\label{erg-eq65}
\lint_{|y| \leq 1} |y|^{\beta} \, \nu(x_0,dy)<\infty \implies \beta_{\infty}^{x_0} \leq \beta.
\end{equation}
In this case, the Blumenthal--Getoor index can be equivalently characterized in terms of fractional moments of the L\'evy measure
\begin{equation*}
\beta_{\infty}^{x_0} = \inf \left\{\gamma>0; \lint_{|y| \leq 1} |y|^{\gamma} \, \nu(x_0,dy)< \infty \right\};
\end{equation*}
this is a special case of \cite[Proposition 5.4]{rs97}, see also \cite{blumen}. \par

For later reference we state the following result which can be found in \cite[Theorem 2.44]{ltp}.

\begin{theorem} \label{def-5}
    Let $(X_t)_{t \geq 0}$ be a rich L\'evy-type process with symbol $q$ and characteristics $(b,Q,\nu)$. Then $(X_t)_{t \geq 0}$ is a semimartingale and its semimartingale characteristics $(B,C,\mu)$ relative to the truncation function $y\I_{(0,1)}(|y|)$ are given by
    \begin{equation}\label{def-eq10}
		B_t = \int_0^t b(X_s) \, ds, \qquad
        C_t = \int_0^t Q(X_s) \, ds, \qquad
        \mu(\cdot,ds,dy) = \nu(X_s,dy) \, ds.
	\end{equation}
\end{theorem}

\section{Main results}

In this section, we present the main results and some illustrating examples. We will point the reader to further results and examples which can be found in Section~\ref{s-erg} and \ref{s-app}.

Our first main result gives a condition on the regularity of a function $f: \mbb{R}^d \to \mbb{R}$ at a fixed point $x_0 \in \mbb{R}^d$ which ensures that the pointwise limit
\begin{equation}\label{main-eq3}
	\lim_{t \to 0} \frac{\mbb{E}^{x_0} f(X_t)-f(x_0)}{t}
\end{equation}
exists. The required regularity is expressed in terms of the generalized Blumenthal--Getoor index $\beta_{\infty}^{x_0}$ of the L\'evy-type process $(X_t)_{t \geq 0}$, cf.\ \eqref{erg-eq4} for the definition.

\begin{theorem}[Regularity at $x_0$]\label{erg-9}
	Let $(X_t)_{t \geq 0}$ be a rich L\'evy-type process with symbol $q$ and characteristics $(b,Q,\nu)$. Suppose that $f \in C_0(\mbb{R}^d)$ satisfies one of the following conditions for some fixed $x_0 \in \mbb{R}^d$.
	\begin{enumerate}[label*=\upshape (A\arabic*),ref=\upshape A\arabic*]
		\item\label{A1}
		There exist constants $\alpha \in [0,2]$, $\alpha>\beta_{\infty}^{x_0}$, and $C>0$ such that
		\begin{equation*}
		  |f(x)-f(x_0)| \leq C |x-x_0|^{\alpha} \fa x \in \ball{x_0}{1}.
		\end{equation*}
		
		\item\label{A2}
		$f$ is differentiable at $x=x_0$ and there exist $\alpha \in [1,2]$, $\alpha>\beta_{\infty}^{x_0}$, and $C>0$ such that
		\begin{equation*}
		|f(x)-f(x_0)-\nabla f(x_0) \cdot (x-x_0)| \leq C |x-x_0|^{\alpha} \fa x \in \ball{x_0}{1}.
		\end{equation*}
		
		\item\label{A3}
		$f$ is twice continuously differentiable in a neighbourhood of $x_0$.
	\end{enumerate}
	Then the limit \begin{equation*}
	\lim_{t \to 0} \frac 1t \left(\mbb{E}^{x_0} f(X_t)-f(x_0)\right)
	\end{equation*}
	exists and takes the value
	\begin{enumerate}
		\item[\eqref{A1}]
		$Lf(x_0) := \lint  \left(f(x_0+y)-f(x_0)\right)  \nu(x_0,dy),$
		\item[\eqref{A2}]
		$Lf(x_0) :=  b(x_0) \cdot \nabla f(x_0) + \lint  \left(f(x_0+y)-f(x_0)-\nabla f(x_0) \cdot y \I_{(0,1)}(|y|)\right)  \nu(x_0,dy),$
		\item[\eqref{A3}]
		$\begin{aligned}[t]\textstyle
		Lf(x_0) &\textstyle:= b(x_0) \cdot \nabla f(x_0) +\tfrac{1}{2}\tr\left(Q(x_0) \cdot \nabla^2 f(x_0)\right)\\
		&\textstyle\qquad\mbox{}+ \lint  \left(f(x_0+y)-f(x_0)-\nabla f(x_0) \cdot y \I_{(0,1)}(|y|)\right)  \nu(x_0,dy),
		\end{aligned}$
	\end{enumerate}
	depending on which of the conditions \eqref{A1}-\eqref{A3} is satisfied.
\end{theorem}

For the proof of Theorem~\ref{erg-9} see Section~\ref{s-erg}, p.~\pageref{proof-erg-9}. As a by-product of the proof, we will find that for any rich L\'evy-type process $(X_t)_{t \geq 0}$ the family of measures $p_t(dy) := t^{-1} \mbb{P}^{x_0}(X_t-x_{0} \in dy)$, $t>0$, on $(\mbb{R}^d \setminus \{0\}, \mc{B}(\mbb{R}^d \setminus \{0\}))$ converges vaguely to $\nu(x_0,dy)$ for each fixed $x_0 \in \mbb{R}^d$, cf.\ Corollary~\ref{erg-7}.

In Theorem~\ref{erg-9} we have to assume that $\alpha$ is strictly larger than the Blumenthal--Getoor index $\beta_{\infty}^{x_0}$. It turns out that Theorem~\ref{erg-9} also holds for $\alpha = \beta_{\infty}^{x_0}$ if $q(x_0,\cdot)$ satisfies a sector condition, has no diffusion part, and the fractional moment $\lint_{|y| \leq 1} |y|^{\beta_{\infty}^{x_0}} \, \nu(x_0,dy)$ is finite, see Theorem~\ref{erg-10} for the precise statement. Moreover, it is possible to extend Theorem~\ref{erg-9} to functions $f$ which are not necessarily bounded, see Theorem~\ref{erg-13} and Theorem~\ref{erg-17} for details.

\medskip
The other main results concern the existence of the limit \eqref{main-eq3} \emph{uniformly} with respect to $x_0 \in \mbb{R}^d$. For the particular case that $(X_t)_{t \geq 0}$ is a L\'evy process we obtain the following statement, see Section~\ref{s-app}, p.~\pageref{proof-app-3}, for the proof.

\begin{theorem} \label{app-3}
	Let $(L_t)_{t \geq 0}$ be a L\'evy process with L\'evy triplet $(b,Q,\nu)$. Denote by $(A,\mc{D}(A))$ its generator and fix $\alpha \in [0,2]$ such that $\lint_{|y| \leq 1} |y|^{\alpha} \, \nu(dy)<\infty$.
	\begin{enumerate}
		\item\label{app-3-i}
		$C_0^2(\mbb{R}^d) \subseteq \mc{D}(A)$ and for $f \in C_0^2(\mbb{R}^d)$
        \begin{equation*}
            Af = b \cdot \nabla f + \tfrac{1}{2} \tr(Q \nabla^2 f) + \lint \left(f(\bullet+y)-f-\nabla f \cdot y \I_{(0,1)}(|y|)\right)  \nu(dy).
        \end{equation*}
		
		\item\label{app-3-ii}
        If $Q=0$, $\alpha \in [0,1]$ and $b=\int_{|y| <1} y \, \nu(dy)$, then $\mc{C}_0^{\alpha}$ is contained in $\mc{D}(A)$ and for $f \in \mc{C}_0^{\alpha}$
		\begin{equation*}
            Af(x) = \lint (f(x+y)-f(x)) \, \nu(dy).
        \end{equation*}
		
		\item\label{app-3-iii}
		If $Q=0$ and $\alpha \in [1,2]$, then $\mc{C}_0^{1,\alpha-1}$ is contained in $\mc{D}(A)$ and for $f \in \mc{C}_0^{1,\alpha-1}$
        \begin{equation*}
			Af(x) = b \cdot \nabla f(x)  + \lint \left(f(x+y)-f(x)-\nabla f(x) \cdot y \I_{(0,1)}(|y|)\right)  \nu(dy).
		\end{equation*}
	\end{enumerate}
\end{theorem}

We refer the reader to \eqref{def-eq2} for the definition of the H\"{o}lder spaces $\mc{C}_0^{\alpha}$ and $\mc{C}_0^{1,\alpha-1}$. Part (ii) of Theorem~\ref{app-3} was recently proved by Cygan \& Grzywny \cite{cygan} for the particular case $\alpha=1$. Let us illustrate Theorem~\ref{app-3} with two examples.

\begin{example}[Isotropic $\alpha$-stable L\'evy processes] \label{app-4}
	Let $(L_t)_{t \geq 0}$ be an isotropic $\alpha$-stable process for some $\alpha \in (0,2)$, i.\,e.\ a L\'evy process with characteristic exponent $\psi(\xi) = |\xi|^{\alpha}$, $\xi \in \mbb{R}^d$, and set $c_{\alpha} :=  \alpha 2^{\alpha-1} \pi^{-d/2} \Gamma\left(\frac{\alpha+d}{2}\right)\big/\Gamma\left(1-\frac\alpha 2\right)$. Then, by Theorem~\ref{app-3}:
	\begin{itemize}
		\item
		If $\alpha \in (0,1)$, then  Theorem~\ref{app-3} shows that the H\"{o}lder space $\mc{C}_0^{\beta}$ is contained in the domain of the generator $A$ for any $\beta \in (\alpha,1]$ and
		\begin{equation*}
		Af(x) = c_{\alpha} \lint \left(f(x+y)-f(x)\right) \frac{dy}{|y|^{d+\alpha}},
		\qquad f \in \mc{C}_0^{\beta},\; x \in \mbb{R}^d.
		\end{equation*}
		\item
		If $\alpha \in [1,2)$, then $\mc{C}_0^{1,\beta-1} \subseteq \mc{D}(A)$ for all $\beta \in (\alpha,2]$ and
		\begin{equation*}
		Af(x) = c_{\alpha} \lint \left(f(x+y)-f(x)-\nabla f(x) \cdot y \I_{(0,1)}(|y|)\right) \frac{dy}{|y|^{d+\alpha}},
		\qquad f \in \mc{C}_0^{1,\beta-1},\; x \in \mbb{R}^d.
		\end{equation*}
	\end{itemize}
The generator of an isotropic $\alpha$-stable L\'evy process is the fractional power $-(-\Delta)^{\alpha/2}$ of the Laplace-operator $\Delta$; this is a well-known fact which probably goes back to Bochner \cite{bochner49} and \cite[p.~93 and pp.~102--106]{bochner55}. Depending on the domain, there are various (equivalent) ways to define fractional powers and we refer to the survey paper \cite{kwas17}. Along with the information on the domain, our example recovers the classical integro-differential representation of the fractional Laplacian as it is widely used in analysis, see e.g.\ \cite{caf-sil14,jacob123,stein} to mention but a few.

Let us mention that the domain $\mc{D}(A)$ of the generator of $(L_t)_{t \geq 0}$ is contained in the Zygmund--H\"{o}lder space $\mathfrak{C}_0^{\alpha} := \mathfrak{C}^{\alpha} \cap C_0$, see Section~\ref{def} for the definition. In dimension $d=1$ this follows by combining two results from interpolation theory \cite[Theorem~1(a), p.~201; Theorem~(d), p.~101]{triebel-inter} with the fact that the domain of the generator of one-dimensional Brownian motion equals $C_0^2(\mbb{R})$ \cite[Example~7.15]{bm2}. For $d \geq 1$ it is possible to show that the resolvent $R_{\lambda}$, $\lambda>0$, satisfies $R_{\lambda}(C_0(\mbb{R}^d)) \subseteq \mathfrak{C}_0^{\alpha}$ using well-known heat kernel estimates for the transition density of $(L_t)_{t \geq 0}$, see e.\,g.\ \cite{blumenthal60} or \cite[formula (2.11)]{kwas17}; since $\mc{D}(A) = R_{\lambda}(C_0(\mbb{R}^d))$ this gives the assertion. In summary,
\begin{equation}
	\mathfrak{C}^{\alpha+}_0 := \bigcup_{\eps>0} \mathfrak{C}_0^{\alpha +\eps}
	\subseteq \mc{D}(A)
	\subseteq \mathfrak{C}_0^{\alpha}. \label{app-eq5}
	\end{equation}
\end{example}

\begin{example}[Compound Poisson processes] \label{app-45}
	Let $(L_t)_{t \geq 0}$ be a L\'evy process with L\'evy triplet $(b,0,\nu)$. Suppose that $\nu$ is a finite measure and $b=\lint_{|y|<1} y \, \nu(dy)$  (e.\,g.\ $b=0$ and $\nu|_{\ball{0}{1}}$ symmetric). Then the domain $\mc{D}(A)$ of the generator of $(L_t)_{t \geq 0}$ equals $C_0(\mbb{R}^d)$ and
	\begin{equation*}
	Af(x) = \lint  \left(f(x+y)-f(x)\right)  \nu(dy),
	\qquad f \in C_0(\mbb{R}^d),\; x \in \mbb{R}^d.
	\end{equation*}
\end{example}

Our third, and final, main result extends Theorem~\ref{app-3} to the much larger class of L\'evy-type processes, see Section~\ref{s-app}, p.~\pageref{proof-app-5} for the proof.

\begin{theorem} \label{app-5}
	Let $(X_t)_{t \geq 0}$ be a rich L\'evy-type process with symbol $q$ and characteristics $(b,Q,\nu)$. Assume that $(X_t)_{t \geq 0}$ has bounded coefficients and that $x \mapsto Q(x)$ is continuous. For fixed $x \in \mbb{R}^d$ denote by $\beta_{\infty}^{x}\in [0,2]$ the generalized  Blumenthal--Getoor index at $\infty$, cf.\ \eqref{erg-eq4}. Let $\alpha: \mbb{R}^d \to (0,2]$ be a uniformly continuous mapping such that $\alpha(x) \geq \min\{\beta_{\infty}^{x}+\eps,2\}$ and
	\begin{equation*}
	\sup_{x \in \mbb{R}^d} \lint_{|y| \leq 1} |y|^{\alpha(x)-\eps} \, \nu(x,dy)<\infty
	\end{equation*}
	for some absolute constant $\eps \in (0,\inf_{x \in \mbb{R}^d} \alpha(x))$. Suppose that $f \in C_0(\mbb{R}^d)$ satisfies the following conditions.
	\begin{enumerate}[label*=\upshape (C\arabic*),ref=\upshape C\arabic*]
		\item\label{C1}
		For any $x \in \{0<\alpha \leq 1\}$ it holds that
		\begin{equation*}
		\sup_{0<|y| \leq 1} \frac{|f(x+y)-f(x)|}{|y|^{\alpha(x)}} < \infty.
		\end{equation*}
		
		\item\label{C2}
		$f$ is differentiable at every point $x \in \{1<\alpha<2\}$ and $g_j(x) := \partial_{x_j} f(x)$, $x \in \{1<\alpha<2\}$, has a $C_0$-extension to $\mbb{R}^d$ for each $j \in \{1,\ldots,d\}$. Moreover,
		\begin{equation*}
		\sup_{0<|y| \leq 1} \frac{|f(x+y)-f(x)-\nabla f(x) \cdot y|}{|y|^{\alpha(x)}} < \infty \fa x \in \{1<\alpha<2\}.
		\end{equation*}
		
		\item\label{C3}
		For any $x \in \{\alpha=2\}$, $f$ is twice differentiable in a neighbourhood of $x$ and the function $h_{ij}(x) := \partial_{x_i} \partial_{x_j} f(x)$, $x \in \{\alpha=2\}$, has a $C_0$-extension to $\mbb{R}^d$ for all $i,j \in \{1,\ldots,d\}$.
	\end{enumerate}
	Then $f$ is in the domain $\mc{D}(A)$ of the generator $A$ of $(X_t)_{t \geq 0}$ and
	\begin{equation*}
	Af(x)
	= b(x) \cdot g(x) + \frac{1}{2} \tr\left(Q(x) h(x)\right) + \lint \left(f(x+y)-f(x)-g(x) \cdot y \I_{(0,1)}(|y|)\right)  \nu(x,dy)
	\end{equation*}
	for all $x \in \mbb{R}^d$ where $g:=(g_1,\ldots,g_d)^\top$ and $h:= (h_{ij})_{i,j=1,\ldots,d}$.
\end{theorem}

As an immediate consequence of Theorem~\ref{app-5} we obtain that certain H\"{o}lder spaces of variable order are contained in the domain of the generator, see Corollary~\ref{app-11} and Corollary~\ref{app-13}. Theorem~\ref{app-5} applies to a large class of L\'evy-type processes, for instance, stable-dominated processes (Example~\ref{app-15}), stable-like processes (Example~\ref{app-17}) and solutions to L\'evy-driven SDEs (Example~\ref{app-19}). In particular, Theorem~\ref{app-5} allows us to obtain the following natural generalization of Example~\ref{app-4}.

\begin{example}[Isotropic stable-like process] \label{main-5}
	Let $(X_t)_{t \geq 0}$ be a rich L\'evy-type process with symbol $q(x,\xi) = |\xi|^{\gamma(x)}$ for a H\"{o}lder continuous mapping $\gamma: \mbb{R}^d \to (0,2)$ which is bounded away from zero. Let $\alpha: \mbb{R}^d \to [0,2]$ be a uniformly continuous mapping such that $\inf_{x \in \mbb{R}^d} (\alpha(x)-\gamma(x))>0$.
	\begin{enumerate}
		\item
		If $\alpha(\mbb{R}^d) \subseteq [0,1]$, then \begin{equation*}
		\mc{C}_0^{\alpha(\cdot)}
		:= \left\{f \in C_0(\mbb{R}^d);\: \sup_{x \in \mbb{R}^d} \sup_{0<|y| \leq 1} \frac{|f(x+y)-f(x)|}{|y|^{\alpha(x)}}<\infty \right\}
		\end{equation*}
		is contained in the domain $\mc{D}(A)$ of the infinitesimal generator $A$ and \begin{equation*}
		Af(x) = c_{\gamma(x)} \lint \left(f(x+y)-f(x)\right) \frac{1}{|y|^{d+\gamma(x)}} \, dy \fa f \in \mc{C}_0^{\alpha(\cdot)}
		\end{equation*}
		where $c_{\gamma(x)} := \gamma(x) \pi^{-d/2} \Gamma((\gamma(x)+d)/2) /\Gamma(1-\gamma(x)/2)$.
		\item The H\"{o}lder space of variable order \begin{equation*}
		\mc{C}_0^{1,(\alpha(\cdot)-1)^+}	
		:=  \left\{f \in C^1_0(\mbb{R}^d); \forall j=1,\ldots,d: \partial_j f \in \mc{C}_0^{\max\{\alpha(\cdot)-1,0\}}\right\}
		\end{equation*}
		is contained in $\mc{D}(A)$ and \begin{equation*}
		Af(x) = c_{\gamma(x)} \lint \left(f(x+y)-f(x)-\nabla f(x) \cdot y \I_{(0,1)}(|y|)\right) \frac{1}{|y|^{d+\gamma(x)}} \, dy  \end{equation*}
		for all $\mc{C}_0^{1,(\alpha(\cdot)-1)^+}$.
	\end{enumerate}
\end{example}

For the existence of the L\'evy-type process $(X_t)_{t \geq 0}$ with symbol $q(x,\xi)=|\xi|^{\gamma(x)}$ we refer the reader to \cite[Theorem 5.2]{matters}, see e.\,g.\ also \cite{bass,kolokoltsov} for related results. Example~\ref{main-5} applies, in particular, in the L\'evy case, i.\,e.\ if $\gamma(x)$ does not depend on $x$, and therefore it generalizes Example~\ref{app-4}.  \par
Let us close this section with some remarks on Theorem~\ref{app-5}.

\begin{remark} \label{app-6}
	\begin{enumerate}
		\item\label{app-6-iv}
		Depending on the local H\"{o}lder index $\alpha(x)$, the generator $Af(x)$, $f\in\mc{D}(A)$, has the following equivalent representations:
		\begin{itemize}
			\item
			$Af(x) = \lint \left(f(x+y)-f(x)\right)  \nu(x,dy)$ for any $x \in \{0<\alpha \leq 1\}$
			
			\item
			$Af(x)= b(x) \cdot \nabla f(x) + \lint \left(f(x+y)-f(x)-\nabla f(x) \cdot y \I_{(0,1)}(|y|)\right)  \nu(x,dy)$ for any $x \in \{1<\alpha<2\}$
			
			\item
			$\begin{aligned}[t]\textstyle
			Af(x) &\textstyle= b(x) \cdot \nabla f(x) +\tfrac{1}{2}\tr\left(Q(x) \cdot \nabla^2 f(x)\right)\\
			&\textstyle\qquad\mbox{}+ \lint  \left(f(x+y)-f(x)-\nabla f(x) \cdot y \I_{(0,1)}(|y|)\right)  \nu(x,dy)
			\end{aligned}$ \\
			for any $x \in \{\alpha=2\}$.
			
		\end{itemize}
		
		\item\label{app-6-iii}
		Since the regularity of the function $f$ may vary from point to point and the triplet is $x$-dependent, Theorem~\ref{app-5} requires stronger assumptions than in the L\'evy case.
		
		\item\label{app-6-ii}
		Let $q$ be a negative definite symbol with characteristics $(b,0,\nu)$ and suppose that $q$ satisfies the sector condition, i.\,e.\ there exists a constant $C>0$ such that
		\begin{equation}\label{app-eq11}
		|\im q(x,\xi)| \leq C \re q(x,\xi) \fa x,\xi \in \mbb{R}^d.
		\end{equation}
		Then $\lint_{|y| \leq 1} |y|^{\alpha(x)-\eps} \, \nu(x,dy)<\infty$ entails $\beta_{\infty}^{x} \leq \alpha(x)-\eps$, cf.\ \eqref{erg-eq65}. Consequently, it suffices in this case to check the integrability condition \begin{equation*}
		\sup_{x \in \mbb{R}^d} \lint_{|y| \leq 1} |y|^{\alpha(x)-\eps} \, \nu(x,dy)<\infty.
		\end{equation*}
		On the other hand, if there exist constants $C>0$ and $\delta>0$ such that
		\begin{equation}\label{app-eq12}
		|\re q(x,\xi)| \leq C |\xi|^{\beta_{\infty}^x+\delta} \fa x,\xi \in \mbb{R}^d,\; |\xi| \geq 1,
		\end{equation}
		then any uniformly continuous function $\alpha: \mbb{R}^d \to (0,2)$ with \begin{equation*}
			\inf_{x \in \mbb{R}^d} (\alpha(x)-\beta_{\infty}^x-\delta)>0
		\end{equation*}
		satisfies the assumptions of Theorem~\ref{app-5}; this follows from the inequality
		\begin{equation*}
		\int_{|y| \leq 1} |y|^{\kappa} \, \nu(dy) \leq c_{\kappa} \int \frac{\re \psi(\xi)}{|\xi|^{d+\kappa}} \, d\xi,
		\qquad \kappa \in (0,2)
		\end{equation*}
		which holds for any continuous negative definite function $\psi: \mbb{R}^d \to \mbb{C}$ with triplet $(b,0,\nu)$, see \eqref{aux-eq1} in the proof of Lemma~\ref{aux-1}. 
		
		The sector condition \eqref{app-eq11} is, in particular, satisfied if $q(x,\cdot)$ is real-valued. This is equivalent to saying that $q(x,\cdot)$ symmetric for all $x \in \mbb{R}^d$ (i.\,e.\ $q(x,\xi) = q(x,-\xi)$ for all $x,\xi \in \mbb{R}^d$) or $b(x)=0$ and $\nu(x,dy)=\nu(x,-dy)$ for all $x\in\rd$.
		
		\item\label{app-6-i}
		It is well known, cf.\ \cite[Theorem 2.30]{ltp}, that the mapping $x \mapsto q(x,\xi)$  is continuous for all $\xi \in \mbb{R}^d$ for any symbol $q$ with $q(x,0)=0$. However, continuity of $q(\cdot,\xi)$ does, in general, not imply continuity of $x \mapsto Q(x)$; consider, for instance, \begin{equation*}
		q(x,\xi) := \frac 12 \xi^2 \I_{\{0\}}(x) + \frac{1-\cos(x \xi)}{x^2} \I_{\rd\setminus\{0\}}(x), \qquad x,\xi \in \mbb{R},
		\end{equation*}
		see \cite[p.\ 11]{courrege}.
	\end{enumerate}
\end{remark}

\section{Pointwise limits} \label{s-erg}

In this section we investigate the small-time asymptotics of generalized moments, i.\,e.\ we study limits of the form
\begin{equation}\label{erg-eq0}
    \lim_{t \to 0} \frac 1t\left(\mbb{E}^x f(X_t)-f(x)\right)
\end{equation}
for a rich Feller process $(X_t)_{t \geq 0}$ and any fixed $x \in \mbb{R}^d$. Recall that a function $f$ is contained in the domain $\mc{D}(A) \subseteq C_0(\mbb{R}^d)$ of the generator $A$, if the limit exists \emph{uniformly} in $C_0(\mbb{R}^d)$, i.\,e.\
\begin{align*}
	\mc{D}(A)
    &:= \left\{ f \in C_0(\mbb{R}^d); \exists g \in C_0(\mbb{R}^d) \,:\,
        \lim_{t \to 0} \sup_{x \in \mbb{R}^d} \left| \frac 1t\left(\mbb{E}^x f(X_t)-f(x)\right) -g(x) \right| = 0 \right\}, \\
	Af(x)
&:= \lim_{t \to 0} \frac 1t\left(\mbb{E}^x f(X_t)-f(x)\right).
\end{align*}
It is, in general, a non-trivial task to check whether a function $f \in C_0(\mbb{R}^d)$ is in the domain of the generator; typically, this requires assumptions on the smoothness, e.\,g.\ $f \in C_0^2(\mbb{R}^d)$ if $(X_t)_{t \geq 0}$ has bounded coefficients, cf.\ \cite[Theorem 2.37(h)]{ltp}.

We are interested in proving the existence of the limit \eqref{erg-eq0} (and also determining it) for functions $f$ which are not necessarily bounded or differentiable. Intuitively, there are two issues which we have to consider:
\begin{enumerate}
    \item We have to ensure that the expectation $\mbb{E}^xf(X_t)$ exists; therefore, we need an assumption on the growth of $f$ at infinity.
    \item For the existence of the limit \eqref{erg-eq0} for a fixed $x \in \mbb{R}^d$ the behaviour of $f$ close to $x \in \mbb{R}^d$ is crucial. For instance, if $X_t := t$ is a deterministic drift process, then the limit \eqref{erg-eq0} exists if, and only if, $f$ is differentiable at $x$. This means that we have to make an assumption on the local regularity of $f$ at $x$, typically H\"{o}lder continuity or differentiability.
\end{enumerate}
In a first step we consider the particular case that $f$ vanishes at infinity and satisfies $f|_{\ball{x}{\delta}}=0$ for some $\delta>0$; for such functions $f$ we show in Theorem~\ref{erg-5}
\begin{equation*}
	\lim_{t \to 0} \frac 1t\left(\mbb{E}^x f(X_t)-f(x)\right) = \lint  \left(f(x+y)-f(x)\right)  \nu(x,dy).
\end{equation*}
This implies, in particular, that $t^{-1} \mbb{P}^x(X_t-x \in \cdot)$ converges vaguely to $\nu(x,\cdot)$ as $t \to 0$, cf.\ Corollary~\ref{erg-7}, and so,
\begin{equation*}
	\lim_{t \to 0} \frac{1}{t} \mbb{P}^x(X_t-x \in A) = \nu(x,A)
\end{equation*}
for any $A \in \mc{B}(\mbb{R}^d \setminus \{0\})$ such that $0\notin\bar A$ and  $\nu(x,\partial A)=0$. In Theorem~\ref{erg-9} and Theorem~\ref{erg-10} we show that the assumption $f|_{\ball{x}{\delta}}=0$ on the regularity of $f$ at $x$ can be replaced by a local H\"{o}lder or differentiability condition. The required regularity can be expressed  in terms of fractional moments of $\nu(x,\cdot)$ or in terms of the generalized Blumenthal--Getoor index at infinity, see \eqref{erg-eq4} for the definition. Finally, in Theorem~\ref{erg-13}, we extend Theorem~\ref{erg-9} to functions $f$ which are not necessarily bounded.

\medskip
The following upper bound for the small-time asymptotics of $\mbb{P}(|X_t-x| \geq r)$ will be one of our main tools.

\begin{lemma} \label{erg-3}
	Let $(X_t)_{t \geq 0}$ be a rich L\'evy-type process with symbol $q$. There exists a function $c: \mbb{R}^d \to (0,\infty)$ such that
    \begin{equation*}
			\limsup_{t \to 0} \frac{1}{t} \mbb{P}^x \left(  |X_t-x| \geq r \right)
			\leq \limsup_{t \to 0} \frac{1}{t} \mbb{P}^x \Big( \sup_{s \leq t} |X_s-x| \geq r \Big)
			\leq c(x) \sup_{|\xi| \leq r^{-1}} |q(x,\xi)|
	\end{equation*}
    for all $x \in \mbb{R}^d$ and $r>0$. Moreover, $c$ is locally bounded, i.\,e.\ $c_K := \sup_{x \in K} c(x)<\infty$ for any compact set $K \subseteq \mbb{R}^d$.
\end{lemma}

Lemma~\ref{erg-3} is a localized variant of a known maximal inequality, cf.\ \cite[Corollary 5.2]{ltp}; for the readers' convenience we include a full proof.

\begin{proof}[Proof of Lemma~\ref{erg-3}]
	For fixed $x \in \mbb{R}^d$ and $r>0$ denote by $\tau_r^x := \inf\{t \geq 0; X_t \notin \ball{x}{r}\}$ the exit time from the ball $\ball{x}{r}$. As
    \begin{equation*}
		\{|X_t-x| \geq r\}
        \subseteq  \Big\{ \sup_{s \leq t} |X_s-x| \geq r \Big\}
        \subseteq \{\tau_r^x \leq t\},
	\end{equation*}
	it suffices to show that
    \begin{equation}\label{erg-st3}
		\limsup_{t \to 0} \frac 1t \mbb{P}^x(\tau_r^x \leq t) \leq c \sup_{|\xi| \leq r^{-1}} |q(x,\xi)| \tag{$\star$}
	\end{equation}
    for some constant $c>0$. To this end, fix $x \in \mbb{R}^d$, $r>0$ and pick $u \in C_c^{\infty}(\mbb{R}^d)$ such that $u(0)=1$, $\spt u \subseteq \ball{0}{1}$ and $0 \leq u \leq 1$. If we set $u_r^x(y) := u((y-x)/r)$, then $u_r^x \in C_c^{\infty}(\mbb{R}^d) \subseteq \mc{D}(A)$, and an application of Dynkin's formula, cf.\ Lemma~\ref{aux-0}, gives
    \begin{equation*}
		\mbb{E}^x u_r^x(X_{t \wedge \tau_r^x}) - 1 =  \mbb{E}^x \left( \int_{[0,t \wedge \tau_r^x)} Au_r^x(X_s) \, ds \right)
	\end{equation*}
	where $A$ denotes the generator of $(X_t)_{t \geq 0}$. Thus, \begin{align*}
		\mbb{P}^x\left(\tau_r^x \leq t\right)
		\leq \mbb{E}^x\left(1-u_r^x(X_{t \wedge \tau_r^x})\right)
		&= -\mbb{E}^x \left( \int_{[0,t \wedge \tau_r^x)} Au_r^x(X_s) \, ds \right) \\
		&= -\mbb{E}^x \left( \int_{[0,t \wedge \tau_r^x)}  \I_{\{|X_{s}-x|<r\}} Au_r^x(X_{s}) \, ds \right).
	\end{align*}
	Since
    \begin{align*}
		-Au_r^x(y)
		= \int_{\mbb{R}^d} e^{iy \cdot \xi} q(y,\xi) \widehat{u_r^x}(\xi) \, d\xi
		&= e^{-ix \cdot \xi} r^d \int_{\mbb{R}^d} e^{iy \cdot \xi} q(y,\xi) \widehat{u}(r \xi) \, d\xi \\
		&= e^{-ix \cdot \xi} \int_{\mbb{R}^d} e^{iy \cdot \xi} q(y,r^{-1} \xi) \widehat{u}(\xi) \, d\xi
	\end{align*}
	for all $y \in \mbb{R}^d$, we get
    \begin{align*}
		\mbb{P}^x\left(\tau_r^x \leq t\right)
		&\leq t \, \mathbb{E}^x \left( \int_{\mbb{R}^d} \sup_{s < t \wedge \tau_r^x} |q(X_{s},r^{-1} \xi)| \, |\widehat{u}(\xi)| \, d\xi \right).
	\end{align*}
    As $X_{s} \in \ball{x}{r}$ for all $s < t \wedge \tau_r^x$, there exists by Theorem~\ref{def-3} a constant $C=C(r,x)$ such that
    \begin{equation*}
		\sup_{s < t \wedge \tau_r^x} |q(X_{s},r^{-1} \xi)| \, |\widehat{u}(\xi)| \leq C (1+|\xi|^2) |\widehat{u}(\xi)| \in L^1(d\xi)
	\end{equation*}
    for all $t \geq 0$. On the other hand, $q(x,0)=0$ implies that $x \mapsto q(x,\xi)$ is continuous for all $\xi \in \mbb{R}^d$, see \cite[Theorem 2.30]{ltp}), and therefore
    \begin{equation*}
		\sup_{s < t \wedge \tau_r^x} |q(X_{s},r^{-1} \xi)| \, |\widehat{u}(\xi)| \xrightarrow[]{t \to 0} |q(x,r^{-1} \xi)| \, |\widehat{u}(\xi)|
    \quad\text{for almost all $\xi \in \mbb{R}^d$.}
	\end{equation*}
	Applying the dominated convergence theorem yields
    \begin{equation*}
		\limsup_{t \to 0} \frac 1t \mbb{P}^x\left(\tau_r^x \leq t\right) \leq \int_{\mbb{R}^d} |q(x,r^{-1} \xi)| \, |\widehat{u}(\xi)| \, d\xi.
	\end{equation*}
	Now \eqref{erg-st3} follows using the estimate from Theorem~\ref{def-3}
    \begin{equation*}
		|q(x, r^{-1}\xi)| \leq 2 \sup_{|\eta| \leq r^{-1}} |q(x,\eta)| (1+|\xi|^2) \fa \xi \in \mbb{R}^d, \;r>0.
    \qedhere
	\end{equation*}
\end{proof}

The following result is well known for the particular case that $(X_t)_{t \geq 0}$ is a L\'evy process, see \cite[Proposition 18.2]{berg} or \cite[Lemma 2.16]{ltp}, its extension to L\'evy-type processes is new. 

\begin{theorem}
\label{erg-5}
    Let $(X_t)_{t \geq 0}$ be a rich L\'evy-type process with symbol $q$ and characteristics $(b,Q,\nu)$. Let $f \in C_0(\mbb{R}^d)$ and suppose that $f|_{\ball{x_0}{\delta}} =0$ for some $x_0 \in \mbb{R}^d$ and $\delta>0$. Then
    \begin{equation*}
		\frac 1t \mbb{E}^x f(X_t) \xrightarrow[]{t \to 0} \lint  f(x+y) \, \nu(x,dy)
	\end{equation*}
	uniformly in a neighbourhood of $x_0$. In particular, $x \mapsto \int f(x+y) \, \nu(x,dy)$ is continuous at $x=x_0$.
\end{theorem}
\begin{proof}
    For fixed $\eps>0$ choose $\chi \in C_c^{\infty}(\mbb{R}^d)$ such that $\|f-\chi\|_{\infty} \leq \eps$. Without loss of generality, we may assume that $\chi|_{\ball{x_0}{\delta}}=0$. Obviously,
    \begin{align*}
		\left| \frac 1t \mbb{E}^x f(X_t) - \lint  f(x+y) \, \nu(x,dy) \right|
		&\leq \left| \frac 1t \mbb{E}^x (f-\chi)(X_t) \right|
        + \lint  \left|f(x+y)-\chi(x+y)\right| \, \nu(x,dy)\\
		&\quad\mbox{}+ \left| \frac 1t \mbb{E}^x \chi(X_t) - \lint  \chi(x+y) \, \nu(x,dy) \right| \\
		&=: I_1+I_2+I_3.
	\end{align*}
    We estimate the terms separately. Using that $\chi(x) =0$, $\nabla \chi(x)=0$ and $\nabla^2 \chi(x)=0$ for all $x \in \ball{x_0}{\delta/4}$, we find for all $x \in \ball{x_0}{\delta/4}$
    \begin{align*}
		I_3
        = \left| \frac1t (\mbb{E}^x \chi(X_t)-\chi(x)) - A\chi(x) \right|
        \leq \sup_{x \in \mbb{R}^d} \left| \frac 1t (\mbb{E}^x\chi(X_t)-\chi(x))-A\chi(x) \right| \xrightarrow[]{t \to 0} 0
	\end{align*}
	as $\chi \in C_c^{\infty}(\mbb{R}^d) \subseteq \mc{D}(A)$. For $I_2$ we note that for any $x \in \ball{x_0}{\delta/4}$ \begin{align*}
		I_2
		\leq \int_{|y| \geq \delta/4} |f(x+y)-\chi(x+y)| \, \nu(x,dy)
    \leq \eps \sup_{x \in \ball{x_0}{\delta/4}} \nu(x,\mbb{R}^d \setminus \ball{0}{\delta/4}).
	\end{align*}
    Note that the constant on the right-hand side is finite, see e.\,g.\ \cite[Theorem 2.30(d)]{ltp}, and $\delta>0$ is a fixed constant which does not depend on $\eps$. Since
    \begin{align*}
		I_1
		\leq \frac{\eps}{t} \mbb{P}^x \left( |X_t-x_0| \geq \frac{\delta}{2} \right)
		\leq \frac{\eps}{t} \mbb{P}^x \left( |X_t-x| \geq \frac{\delta}{4} \right)
	\end{align*}
	for all $x \in \ball{x_0}{\delta/4}$, it follows from Lemma~\ref{erg-3} that there exists a constant $C>0$ such that
	\begin{equation*}
		\limsup_{t \to 0} \frac{1}{t} I_1
		\leq C  \eps \sup_{x \in \ball{x_0}{\delta/4}} \sup_{|\xi| \leq 4 \delta^{-1}} |q(x,\xi)|.
	\end{equation*}
	The above estimates show
    \begin{align*}
		\limsup_{t \to 0} \bigg| \frac 1t \mbb{E}^x &f(X_t) - \lint  f(x+y) \, \nu(x,dy) \bigg| \\
		&\leq \eps \Big(\sup_{x \in \ball{x_0}{\delta/4}} \nu(x,\mbb{R}^d \setminus \ball{0}{\delta/4})
            + C \sup_{x \in \ball{x_0}{\delta/4}} \sup_{|\xi| \leq 2 \delta^{-1}} |q(x,\xi)| \Big)
		\xrightarrow[]{\eps \to 0} 0.
	\end{align*}
    The assertion on the continuity follows directly from the local uniform convergence and the fact that $x \mapsto \mbb{E}^x f(X_t)$ is continuous as $(X_t)_{t \geq 0}$ is a Feller process.
\end{proof}

If we use Theorem~\ref{erg-5} for the shifted function $f(\cdot -x_0)$ for a fixed $x_0 \in \mbb{R}^d$, we get:

\begin{corollary} \label{erg-7}
    Let $(X_t)_{t \geq 0}$ be a rich L\'evy-type process with symbol $q$ and characteristics $(b,Q,\nu)$. If $f \in C_0(\mbb{R}^d)$ and $f|_{\ball{0}{\delta}} = 0$ for some $\delta>0$, then
    \begin{equation*}
		\lim_{t \to 0} \frac 1t \mbb{E}^x f(X_t-x) = \lint f(y) \, \nu(x,dy) \fa x \in \mbb{R}^d.
	\end{equation*}
\end{corollary}

Corollary~\ref{erg-7} shows that the family of measures $p_t(dy) := t^{-1} \mbb{P}^x(X_t-x \in dy)$, $t>0$, on $(\mbb{R}^d \setminus \{0\}, \mc{B}(\mbb{R}^d \setminus \{0\}))$ converges vaguely to $\nu(x,dy)$ for each fixed $x \in \mbb{R}^d$. By the portmanteau theorem, Corollary~\ref{erg-7} implies \begin{equation}\label{erg-eq5}
	\lim_{t \to 0} \frac 1t \mbb{P}^x(X_t-x \in A) = \nu(x,A)
\end{equation}
for any Borel set $A \in \mc{B}(\mbb{R}^d \setminus \{0\})$ such that $0\notin\bar A$ and $\nu(x,\partial A)=0$.

We are now ready to prove our first main result, Theorem~\ref{erg-9}, see p.~\pageref{erg-9} for the statement. It allows us to to relax the assumption ``$f|_{\ball{x_0}{\delta}}=0$'' in Theorem~\ref{erg-5}.

\begin{proof}[Proof of Theorem~\ref{erg-9}]
\phantomsection\label{proof-erg-9}Pick $\chi \in C_c^{\infty}(\mbb{R}^d)$, $0 \leq \chi \leq 1$, such that $\chi|_{\ball{x_0}{1}}=1$, $\chi|_{\coball{x_0}{2}}=0$ and set $\chi_{\delta}(x) := \chi(\delta^{-1}x)$ for $\delta>0$. \par
\medskip\noindent\eqref{A1}
Without loss of generality, we may assume $f(x_0)=0$, otherwise we consider the shifted function $x \mapsto f(x)-f(x_0)$. As $\alpha>\beta_{\infty}^{x_0}$, we have 
\begin{equation*}
	\lint  |f(x_0+y)| \, \nu(x_0,dy)
	\leq C\lint_{|y| \leq 1} |y|^{\alpha} \, \nu(x_0,dy) + \|f\|_{\infty} \nu(x_0,\mbb{R}^d \setminus \ball{0}{1}) < \infty,
\end{equation*}
and therefore it follows from Theorem~\ref{erg-5} and the dominated convergence theorem that
\begin{align*}
		\frac 1t \mbb{E}^{x_0}([f (1-\chi_{\delta})](X_t))
		&\xrightarrow[\phantom{\delta \to 0}]{t \to 0} \lint f(x_0+y) (1-\chi_{\delta}(x_0+y)) \, \nu(x_0,dy) \\
		&\xrightarrow[\phantom{t \to 0}]{\delta \to 0} \lint f(x_0+y) \, \nu(x_0,dy).
\end{align*}
On the other hand, if we set $C_{\delta} := \sup_{|y-x_0| \leq 2 \delta} |f(y)|$, then $C_{\delta}\to 0$ as $\delta \to 0$ and
\begin{align*}
		\left|\mbb{E}^{x_0} ([f \chi_{\delta}](X_t))\right|
		&\leq \int_0^{C_{\delta}} \mbb{P}^{x_0}\big(|f(X_t)| \geq r,\: |X_t-x_0| \leq 2\delta\big) \, dr \\
		&\leq \int_0^{C_{\delta}} \mbb{P}^{x_0}\big(|X_t-x_0|^{\alpha} \geq r/C\big) \, dr
\end{align*}
for any $\delta \in (0,1/2)$. By Lemma~\ref{erg-3}
\begin{equation}\label{erg-eq7}\begin{aligned}
		\limsup_{t \to 0} \frac{1}{t} \mbb{P}^{x_0}\big(|X_t-x_0|^{\alpha} \geq r/C\big)
		&= \limsup_{t \to 0} \frac{1}{t} \mbb{P}^{x_0}\big(|X_t-x_0| \geq C^{-1/\alpha} r^{1/\alpha}\big) \\
		&\leq c \sup_{|\xi| \leq r^{-1/\alpha} C^{1/\alpha}} |q(x_0,\xi)|
		\leq C' r^{-\beta/\alpha}
\end{aligned}\end{equation}
for any $\beta \in (\beta_{\infty}^x,\alpha)$ and suitable constants $c,C'>0$; thus, by Fatou's lemma, 
\begin{align*}
    \limsup_{t \to 0} \left| \frac 1t \mbb{E}^{x_0}([f \chi_{\delta}](X_t)) \right|
		\leq C' \int_0^{C_{\delta}} r^{-\beta/\alpha} \,dr
		\xrightarrow[]{\delta \to 0} 0.
\end{align*}
Writing
\begin{equation*}
    \frac 1t \mbb{E}^{x_0} f(X_t) =  \frac 1t \mbb{E}^{x_0}([f \chi_{\delta}](X_t)) + \frac 1t \mbb{E}^{x_0}([f (1-\chi_{\delta})](X_t))
\end{equation*}
and letting first $t \to 0$ and then $\delta \to 0$, proves the claim.

\medskip\noindent\eqref{A2}
For fixed $R>0$ let $\tau_R^{x_0}$ denote the exit time from the ball $\ball{x_0}{R}$. The function
\begin{equation*}
    x \mapsto g(x) := f(x)-f(x_0)-\nabla f(x_0)\cdot (x-x_0) \chi(x)
\end{equation*}
satisfies \eqref{A1} and, therefore, by the first part of this proof,
\begin{align*}
		\lim_{t \to 0} \frac 1t \mbb{E}^{x_0} g(X_t)
		&= \lint  \left(g(x_0+y)-g(x_0)\right)  \nu(x_0,dy) \\
		&= \lint  \left(f(x_0+y)-f(x_0)-\chi(y+x_0) \nabla f(x_0) \cdot y\right)  \nu(x_0,dy).
\end{align*}
As $(\bullet-x_0) \chi(\bullet) \in C_c^{\infty}(\mbb{R}^d) \subseteq \mc{D}(A)$ an application of Dynkin's formula, cf.\ Lemma~\ref{aux-0}, shows
\begin{equation*} 
		\frac{1}{t} \mbb{E}^{x_0}\left((X_{t \wedge \tau_R^{x_0}}-x_0) \chi(X_{t \wedge \tau_R^{x_0}})\right)
		\xrightarrow[]{t \to 0} b(x_0) +  \lint  y \left(\chi(y+x_0)-\I_{(0,1)}(|y|)\right)  \nu(x_0,dy)
\end{equation*}
for any $R>0$. Using the fact that $\spt \chi \subseteq \cball{x_0}{2}$ and applying Lemma~\ref{erg-3}, we find for some constant $c=c(x_0)$
\begin{align*}
		&\left| \frac{1}{t} \mbb{E}^{x_0}\left((X_{t \wedge \tau_R^{x_0}}-x_0) \chi(X_{t \wedge \tau_R^{x_0}})\right)
        - \frac{1}{t} \mbb{E}^{x_0}\left((X_t-x_0) \chi(X_t)\right) \right|\\
		&\qquad\leq  \frac 4t \mbb{P}^{x_0}\left(\tau_R^{x_0} \leq t\right)
		\leq 4c \sup_{|\xi| \leq R^{-1}} |q(x_0,\xi)|
		\xrightarrow[]{R \to \infty} 0,
\end{align*}
and therefore we conclude
\begin{equation*}
		\frac 1t \mbb{E}^{x_0} \left((X_t-x_0) \chi(X_t)\right)
        \xrightarrow[]{t \to 0}
        b(x_0)  + \lint  y \left(\chi(y+x_0)-\I_{(0,1)}(|y|)\right)  \nu(x_0,dy).
\end{equation*}
Consequently,
\begin{align*}
		&\frac 1t \left(\mbb{E}^{x_0}f(X_t)-f(x_0)\right)
		= \frac 1t \mbb{E}^{x_0} g(X_t) + \frac 1t \nabla f(x_0) \cdot \mbb{E}^{x_0}\left((X_t-x_0) \chi(X_t)\right) \\
		&\qquad\qquad\xrightarrow[]{t \to 0} b(x_0) \cdot \nabla f(x_0)
        + \lint  \left(f(x_0+y)-f(x_0)-\nabla f(x_0) \cdot y \I_{(0,1)}(|y|)\right)  \nu(x_0,dy),
\end{align*}
finishing the second part.

\medskip\noindent\eqref{A3}
We begin with the particular case that $f(x_0) = 0$ and $\nabla f(x_0)=0$. Since, by Theorem~\ref{erg-5} and the dominated convergence theorem,
\begin{align*}
		\frac{1}{t} \mbb{E}^{x_0} ([f(1-\chi_{\delta})](X_t))
		&\xrightarrow[\phantom{\delta \to 0}]{t \to 0} \lint  [f (1-\chi_{\delta})](x_0+y) \, \nu(x_0,dy)\\
		&\xrightarrow[\phantom{t \to 0}]{\delta \to 0} \lint  f(x_0+y) \, \nu(x_0,dy),
\end{align*}
it is enough to show
\begin{equation}\label{erg-eq9}
    \frac 1t \mbb{E}^{x_0}\left([f \chi_{\delta}](X_t)\right)
    \xrightarrow[]{t,\delta \to 0} \sum_{i,j=1}^d Q_{ij}(x_0) \partial_i \partial_j f(x_0).
\end{equation}
In order to keep notation simple, we set $f_{\delta}(x) := f(x) \chi_{\delta}(x)$. Note that by Lemma~\ref{erg-3}
\begin{align*}
		\left| \frac 1t \mbb{E}^{x_0}f_{\delta}(X_t) - \frac 1t \mbb{E}^{x_0} f_{\delta}(X_{t \wedge \tau_R^{x_0}}) \right|
		&\leq 2 \|f\|_{\infty} \frac 1t \mbb{P}^{x_0}\left(\tau_R^{x_0} \leq t\right) \\
		&\leq 2c \|f\|_{\infty} \sup_{|\xi| \leq R^{-1}} |q(x_0,\xi)|
		\xrightarrow[]{R \to \infty} 0,
\end{align*}
and therefore \eqref{erg-eq9} follows if we can show that
\begin{equation}\label{erg-eq11}
    \frac 1t \mbb{E}^{x_0}(f_{\delta}(X_{t \wedge \tau_R^x}))
    \xrightarrow[]{t,\delta \to 0} \sum_{i,j=1}^d Q_{ij}(x_0) \partial_i \partial_j f(x_0)
\end{equation}
for every fixed $R>0$. By Taylor's formula, there exists a continuous mapping $\varphi: \mbb{R} \to \mbb{R}$ such that $\lim_{r \to 0} \varphi(r)=0$ and
\begin{equation*}
		f(y) = \frac{1}{2} \sum_{i,j=1}^d (y^i-x_0^i)(y^j-x_0^j) \partial_i \partial_j f(x_0) + |y-x_0|^2 \varphi(|x_0-y|)
\end{equation*}
for all $y=(y^1,\ldots,y^d) \in \ball{x_0}{\delta}$. Thus,
\begin{equation*}
    \frac 1t \mbb{E}^{x_0}(f_{\delta}(X_{t \wedge \tau_R^x}))
    = I_1+I_2
\end{equation*}
where
\begin{align*}
		I_1
        &:= \frac{1}{2t} \sum_{i,j=1}^d \partial_i \partial_j f(x_0) \mbb{E}^{x_0}\big[(X_{t  \wedge \tau_R^x}^i-x_0^i)(X_{t \wedge \tau_R^x}^j-x_0^j) \chi_{\delta}(X_{t \wedge \tau_R^x})\big]
    \\
		I_2
        &:= \frac{1}{t} \mbb{E}^{x_0} \left[|X_{t \wedge \tau_R^x}-x_0|^2 \varphi(|X_{t \wedge \tau_R^x}-x_0|) \chi_{\delta}(X_{t \wedge \tau_R^x})\right].
\end{align*}
We estimate the terms separately. By the definition of $\chi_{\delta}$, we have
\begin{equation*}
		I_2 \leq t^{-1} \sup_{r \leq 2\delta} |\varphi(r)|\,\mbb{E}^{x_0}(|X_{t \wedge \tau_R^x}-x_0|^2 \chi(X_{t \wedge \tau_R^x})),
\end{equation*}
and so an application of Dynkin's formula yields
\begin{align*}
		I_2
		\leq \sup_{r \leq 2\delta} |\varphi(r)| \sup_{|y-x_0| \leq R} |A(|\bullet-x_0|^2 \cdot \chi(\bullet))(y)|
		\xrightarrow[]{\delta \to 0} 0.
\end{align*}
Using that $\nabla \chi_{\delta}(x_0) = 0$ and $\nabla^2 \chi_{\delta}(x_0)=0$, it is not difficult to see from Dynkin's formula and the fundamental theorem of calculus that
\begin{align*}
	I_1
    &\xrightarrow[\phantom{\delta \to 0}]{t \to 0} \frac{1}{2} \sum_{i,j=1}^d \partial_i \partial_j f(x_0) \left( Q_{ij}(x_0)
        + \lint  y_i y_j \chi_{\delta}(x_0+y) \, \nu(x_0,dy) \right) \\
	&\xrightarrow[\phantom{t \to 0}]{\delta \to 0}  \frac{1}{2} \sum_{i,j=1}^d \partial_i \partial_j f(x_0) Q_{ij}(x_0).
\end{align*}
Combining both convergence results proves \eqref{erg-eq11} if $f(x_0) = 0$ and $\nabla f(x_0)=0$. For the general case  define
\begin{equation*}
		g(x) := f(x)-f(x_0)-\chi(x) \nabla f(x_0) \cdot (x-x_0), \qquad x \in \mbb{R}^d,
\end{equation*}
and use exactly the same reasoning as in the proof of \eqref{A2}. \qedhere
\end{proof}
	
In Theorem~\ref{erg-9} we have to assume that $\alpha$ is strictly larger than the Blumenthal--Getoor index $\beta_{\infty}^{x_0}$ defined in \eqref{erg-eq4}. In fact, Theorem~\ref{erg-9} also holds for $\alpha = \beta_{\infty}^{x_0}$ if $q(x_0,\cdot)$ satisfies the sector condition, has no diffusion part, and the fractional moment $\lint_{|y| \leq 1} |y|^{\beta_{\infty}^{x_0}} \, \nu(x_0,dy)$ is finite.

\begin{theorem} \label{erg-10}
    Let $(X_t)_{t \geq 0}$ be a rich L\'evy-type process with symbol $q$ and characteristics $(b,0,\nu)$. Suppose that $f \in C_0(\mbb{R}^d)$ satisfies one of the following conditions for some fixed $x_0 \in \mbb{R}^d$.
    \begin{enumerate}[label*=\upshape (B\arabic*),ref=\upshape B\arabic*]
        \item\label{B1}
            There exist $\alpha \in (0,1]$ and $C>0$ such that $\lint_{|y| \leq 1} |y|^{\alpha} \, \nu(x_0,dy)<\infty$ and
            \begin{equation*}
        			|f(x)-f(x_0)| \leq C |x-x_0|^{\alpha} \fa x \in \ball{x_0}{1}.
            \end{equation*}

        \item\label{B2}
            $f$ is differentiable at $x=x_0$ and there exist constants $\alpha \in (1,2)$ and $C>0$ such that $\lint_{|y| \leq 1} |y|^{\alpha} \, \nu(x_0,dy)<\infty$ and
            \begin{equation*}
        			|f(x)-f(x_0)-\nabla f(x_0) \cdot (x-x_0)| \leq C |x-x_0|^{\alpha} \fa x \in \ball{x_0}{1}.
        	\end{equation*}
    \end{enumerate}
    If $q(x_0,\cdot)$ satisfies the sector condition, i.\,e.\ $|\im q(x_0,\xi)| \leq C' \re q(x_0,\xi)$ for some constant $C'>0$, then the limit
    \begin{equation*}
    		\lim_{t \to 0} \frac 1t \left(\mbb{E}^{x_0} f(X_t)-f(x_0)\right)
    \end{equation*}
    exists and takes the value
    \begin{enumerate}
        \item[\eqref{B1}]
            $Lf(x_0) := \lint  \left(f(x_0+y)-f(x_0)\right)  \nu(x_0,dy);$

        \item[\eqref{B2}]
            $Lf(x_0) :=  b(x_0) \cdot \nabla f(x_0) + \lint  \left(f(x_0+y)-f(x_0)-\nabla f(x_0) \cdot y \I_{(0,1)}(|y|)\right)  \nu(x_0,dy).$
    \end{enumerate}
\end{theorem}
\begin{proof}
    The proof is very similar to that of Theorem~\ref{erg-9}; the only modification is needed in \eqref{erg-eq7} where we use the fact that $\int_{|y| \leq 1} |y|^{\alpha} \, \nu(x_0,dy)<\infty$ implies
    \begin{equation*}
        \int_{0}^1 \sup_{|\xi| \leq r^{-1/\alpha}} |q(x_0,\xi)| \, dr
        = \alpha \int_1^{\infty} \frac{1}{s^{1+\alpha}} \sup_{|\xi| \leq s} |q(x_0,\xi)| \, ds < \infty
	\end{equation*}
	(cf.\ Lemma~\ref{aux-1} for details) to obtain an integrable majorant.
\end{proof}

In the remaining part of this section we extend Theorem~\ref{erg-9} and Theorem~\ref{erg-10} to functions $f$ which are not necessarily bounded. Recall that a function $g \geq 0$ is called submultiplicative if there exists a constant $c>0$ such that 
$
	g(x+y) \leq c g(x) g(y)
$
holds for all $x,y \in \mbb{R}^d$.
In \cite{kuehn} it was shown that the implication
\begin{equation*}
    \sup_{x \in K} \lint_{|y| \geq 1} g(y) \, \nu(x,dy)< \infty
    \implies \forall t>0 \,:\,\, \sup_{x \in K} \sup_{s \leq t} \mbb{E}^x g(X_{s \wedge \tau_K}-x)<\infty
\end{equation*}
holds for any twice differentiable submultiplicative function $g \geq 0$, any compact set $K \subseteq \mbb{R}^d$, and any rich L\'evy-type process; if $(X_t)_{t \geq 0}$ has bounded coefficients, then $K=\mbb{R}^d$ is admissible.  Here $\tau_K$ denotes as usual the first exit time from $K$. It is therefore a natural idea to replace
\begin{equation*}
	\frac 1t\left(\mbb{E}^x f(X_t)-f(x)\right) \quad \text{by} \quad \frac 1t\left(\mbb{E}^x f(X_{t \wedge \tau_K})-f(x)\right),
\end{equation*}
and to consider functions $f: \mbb{R}^d \to \mbb{R}$ which can be dominated by a submultiplicative function $g \geq 0$ with $\sup_{x \in K} \int_{|y| \geq 1} g(y) \, \nu(x,dy)<\infty$.
		
\begin{definition} \label{erg-11}
    Let $(b(x),Q(x),\nu(x,dy))$ be an $x$-dependent L\'evy triplet and $K \subseteq \mbb{R}^d$. We write $\Sigma(K)$ for the family of twice differentiable submultiplicative functions $g: \mbb{R}^d \to (0,\infty)$ satisfying the following two integrability conditions.
    \begin{enumerate}
		\item $M(K) := \sup_{x \in K} \int_{|y| \geq 1} g(y) \, \nu(x,dy)< \infty$ \qquad\quad (integrability).
		\item $M_R(K) := \sup_{x \in K} \int_{|y| \geq R} g(y) \, \nu(x,dy) \xrightarrow[]{R \to \infty} 0$ \qquad (tightness).
	\end{enumerate}
\end{definition}

\begin{theorem}[Behaviour at $\infty$] \label{erg-13}
    Let $(X_t)_{t \geq 0}$ be a rich L\'evy-type process with symbol $q$ and characteristics $(b,Q,\nu)$. Moreover, let $f: \mbb{R}^d \to \mbb{R}$ be a continuous mapping satisfying the following growth condition \eqref{G}.
    \begin{enumerate}[label*=\upshape (G),ref=\upshape G]
	\item\label{G}
        There exist a compact set $K \subseteq \mbb{R}^d$ and a function $g \in \Sigma(K)$ such that
        \begin{equation*}
		  \lim_{|x| \to \infty} \left| \frac{f(x)}{g(x)} \right| < \infty.
	   \end{equation*}
	\end{enumerate}
	If one of the conditions \eqref{A1}-\eqref{A3} holds for some $x_0 \in K$, then the limit
    \begin{equation*}
		\lim_{t \to 0} \frac 1t \left(\mbb{E}^{x_0} f(X_{t \wedge \tau_K})-f(x_0)\right)
	\end{equation*}
    exists and equals $Lf(x_0)$ defined in Theorem~\ref{erg-9}; here $\tau_K := \inf\{t \geq 0; X_t \notin K\}$ denotes the exit time from the set $K$. If $(X_t)_{t \geq 0}$ has bounded coefficients, then $K=\mbb{R}^d$ is admissible.
\end{theorem}

\begin{proof}
    We only consider the case that $(X_t)_{t \geq 0}$ has bounded coefficients and $g \in \Sigma(\mbb{R}^d)$; the proof of the other assertion works analogously and just requires an additional stopping argument. For simplicity of notation we assume that $b(x)=0$ and $Q(x)=0$ for all $x \in \mbb{R}^d$, see the remark at the end of the proof.

    Let $\chi$ be a continuous function such that $1-\chi \in C_c^{\infty}(\mbb{R}^d)$, $0 \leq \chi \leq 1$, $\chi|_{\ball{0}{1}}=0$ and $\chi|_{\coball{0}{2}}=1$, and set $\chi_R(x) := \chi(R^{-1}x)$. Then $f(\bullet) \cdot (1-\chi_R(\bullet-x_0))$ satisfies the assumptions of Theorem~\ref{erg-9} for each $R>0$ and therefore
    \begin{equation*}
		\frac{1}{t} \left( \mbb{E}^{x_0}(f(X_t) (1-\chi_R)(X_t-x_0))-f(x_0) \right)
        \xrightarrow[]{t \to 0} L(f (1-\chi_R)(\bullet-x_0))(x_0).
	\end{equation*}
    Since $\nabla \chi_R(x_0)=0$, $\nabla^2 \chi_R(x_0)=0$ for each $R>0$ and $\int_{|y| \geq 1} |f(y)| \, \nu(x_0,dy)<\infty$, it follows easily from the definition of $L(f (1-\chi_R)(\bullet-x_0))$ and the dominated convergence theorem that
    \begin{align*}
		\frac{1}{t} \left( \mbb{E}^{x_0}(f(X_t) (1-\chi_R)(X_t-x_0))-f(x_0) \right)
		&\xrightarrow[\phantom{R \to \infty}]{t \to 0} L(f (1-\chi_R)(\bullet-x_0))(x_0) \\
		&\xrightarrow[\phantom{t \to 0}]{R \to \infty} Lf(x_0).
	\end{align*}
	Consequently, it remains to show that
    \begin{equation*}
		\limsup_{R \to \infty} \limsup_{t \to 0} \left| \frac{1}{t} \mbb{E}^{x_0}\left(f(X_t)\chi_R(X_t-x_0)\right) \right| = 0.
	\end{equation*}
	Because of the growth condition \eqref{G} and the submultiplicativity of $g$, it suffices to prove
	\begin{equation}\label{erg-eq21}
		\limsup_{R \to \infty} \limsup_{t \to 0}  \frac{1}{t} \mbb{E}^{x_0}\left(g(X_t-x_0)\chi_R(X_t-x_0)\right) = 0.
	\end{equation}
    By Theorem~\ref{def-5}, $(X_t)_{t \geq 0}$ is a semimartingale with semimartingale characteristics $(0,0,\mu)$ given by \eqref{def-eq10}. Consequently, $(X_t)_{t \geq 0}$ has a canonical representation $X_t = x_0 + X_t^{(1)}+X_t^{(2)}$,
    \begin{align*}
		X_t^{(1)} &:= \int_0^t \!\! \int_{0<|y| < 1}y \, (N(dy,ds)-\mu(dy,ds)) \\
		X_t^{(2)} &:= \int_0^t \!\! \int_{|y| \geq 1} y \, N(dy,ds)
	\end{align*}
    where $N$ denotes the jump measure of $(X_t)_{t \geq 0}$, cf.\ \cite[Theorem II.2.34]{js}. By the submultiplicativity of $g$, there exists a constant $c>0$ such that
    \begin{equation*}
		g(X_t-x_0)= g(X_t^{(1)}+X_t^{(2)}) \leq c g(X_t^{(1)}) g(X_t^{(2)}), \qquad t \geq 0.
	\end{equation*}
	Since any submultiplicative function grows at most exponentially, cf.\ \cite[Lemma 25.5]{sato}, we can find constants $a,b>0$ such that
    \begin{equation}\label{erg-eq22}
		g(X_t-x_0) \leq a \exp \left( b \sqrt{|X_t^{(1)}|^2+1}-1 \right) g(X_t^{(1)}), \qquad t \geq 0.
	\end{equation}
    In order to keep our notation simple, we assume that $a = b = c = 1$. Moreover, we set
    \begin{equation*}
		\varrho(x) := \exp\left(\sqrt{|x|^2+1}-1\right)
	\end{equation*}
	and use the subscript to denote truncated functions, e.\,g.\
    \begin{equation*}
		\varrho_R(x) := \chi_R(x) \varrho(x) \quad \text{and} \quad g_R(x) := \chi_R(x) g(x).
	\end{equation*}
	From the definition of $\chi_R$ and the triangle inequality, it is not difficult to see that
    \begin{equation}\label{erg-eq24}
		\chi_R(x+y) \leq \chi_{R/4}(x)+\chi_{R/4}(y) \fa x,y \in \mbb{R}^d,
	\end{equation}
	and therefore we obtain
    \begin{align*}
		g(X_t-x_0)\chi_R(X_t-x_0)
		&\leq \exp \left( \sqrt{|X_t^{(1)}|^2+1}-1\right) g(X_t^{(2)})\chi_{R/4}(X_t^{(1)}) \\
		&\quad + \exp \left(\sqrt{|X_t^{(1)}|^2+1}-1 \right) g(X_t^{(2)})\chi_{R/4}(X_t^{(2)}) \\
		&= \varrho_{R/4}(X_t^{(1)}) g(X_t^{(2)}) + \varrho(X_t^{(1)}) g_{R/4}(X_t^{(2)}).
	\end{align*}
	Consequently, \eqref{erg-eq21} follows if we can show \begin{align}
		\lim_{R \to \infty} \lim_{t \to 0} \frac{1}{t} \mbb{E}^{x_0}\left(\varrho_{R/4}(X_t^{(1)}) g(X_t^{(2)})\right) &= 0 \label{erg-eq23} \\
		\lim_{R \to \infty} \lim_{t \to 0} \frac{1}{t} \mbb{E}^{x_0}\left(\varrho(X_t^{(1)}) g_{R/4}(X_t^{(2)})\right) &= 0 \label{erg-eq25}.
	\end{align}
	First we prove \eqref{erg-eq23}. Define a stopping time by
    \begin{equation*}
		\tau := \tau_r := \inf\left\{t>0; |X_t^{(1)}| + |X_t^{(2)}| \geq r\right\}
	\end{equation*}
	for fixed $r>0$. Applying It\^o's formula for semimartingales gives
    \begin{align}
		&\mbb{E}^{x_0}\big(\varrho_{R/4}(X_{t \wedge \tau}^{(1)})  g(X_{t \wedge \tau}^{(2)})\big) \notag \\
        &= \mbb{E}^{x_0} \left( \int_0^{t \wedge \tau} \!\! \int_{|y| \geq 1} \varrho_{R/4}(X_{s}^{(1)}) (g(X_s^{(2)}+y)-g(X_s^{(2)})) \, \nu(X_{s},dy) \, ds \right) \label{erg-st5}
    \\
        &\,\,\mbox{} + \mbb{E}^{x_0} \left( \int_0^{t \wedge \tau} \!\! \lint_{|y|<1} g(X_s^{(2)}) \left(\varrho_{R/4}(X_s^{(1)}+y)-\varrho_{R/4}(X_s^{(1)})-\nabla \varrho_{R/4}(X_s^{(1)}) \cdot y\right)  \nu(X_{s},dy) \, ds \right). \notag
	\end{align}
    Since $g \geq 0$ is submultiplicative, the first term on the right-hand side of \eqref{erg-st5} is bounded above by
    \begin{align*}
		\mbb{E}^{x_0} &\left( \int_0^t \!\! \lint \varrho_{R/4}(X_s^{(1)}) g(X_s^{(2)}) g(y) \, \nu(X_{s},dy) \, ds \right)\\
        &\leq \left(\sup_{x \in \mbb{R}^d} \int_{|y| \geq 1} |g(y)| \, \nu(x,dy)\right) \mbb{E}^{x_0} \left( \int_0^t \varrho_{R/4}(X_s^{(1)}) g(X_s^{(2)}) \, ds \right).
	\end{align*}
    For the second term in \eqref{erg-st5} we apply Taylor's formula and use the fact that $\nabla^2 \chi_{R/4}(z) =0$ for all $z \in \ball{0}{R/4} \cup \coball{0}{R/2}$ to conclude that there exists a function $\psi \in C_b^2(\mbb{R}^d)$ such that $\psi(z)=0$ for all $z \in \ball{0}{1/16}$ and \begin{equation*}
		|\varrho_{R/4}(x+y)-\varrho_{R/4}(x)-\nabla \varrho_{R/4}(x) \cdot y|
        \leq |y|^2 \varrho(x) \psi(x)
        \fa x \in \mbb{R}^d, \; |y| \leq 1
	\end{equation*}
    for $R \geq 1$. Using this estimate for $x := X_{s}^{(1)}$, we find that the second term on the right-hand side of \eqref{erg-st5} is bounded above by
    \begin{align*}
        \left( \sup_{x \in \mbb{R}^d} \lint_{|y| \leq 1} |y|^2 \, \nu(x,dy) \right) \mbb{E}^{x_0} \left( \int_0^t \varrho(X_s^{(1)}) \psi(X_s^{(1)}) g(X_s^{(2)}) \, ds \right).
	\end{align*}
	Now it follows from Fatou's lemma, Definition~\ref{erg-11} and Lemma~\ref{erg-15} below that there exists an absolute constant $C>0$ such that
    \begin{align*}
		\frac{1}{t} \mbb{E}^{x_0}\left(\varrho_{R/4}(X_t^{(1)}) g(X_t^{(2)})\right)
        &\leq \liminf_{r \to \infty} \frac{1}{t} \mbb{E}^{x_0}\left(\varrho_{R/4}(X_{t \wedge \tau}^{(1)})  g(X_{t \wedge \tau}^{(2)})\right)\
		\leq  \frac{C}{t} \int_0^t s \, ds
	\end{align*}
	(recall the definition of $\varrho$, $\varrho_{R/4}$ and note that $K=\mbb{R}^d$), and this implies \eqref{erg-eq23}.

    It remains to prove \eqref{erg-eq25}. Again an application of It\^o's formula shows
    \begin{align}
    	&\mbb{E}^{x_0}(\varrho(X_{t \wedge \tau}^{(1)}) g_{R/4}(X_{t \wedge \tau}^{(2)})) \notag \\
        &=\mbb{E}^{x_0} \left( \int_0^{t \wedge \tau}\!\! \int_{|y| \geq 1} \varrho(X_s^{(1)}) \left(g_{R/4}(X_s^{(2)}+y)-g_{R/4}(X_s^{(2)})\right)  \nu(X_{s-},dy) \, ds \right) \label{erg-eq27} \\
        &\quad\mbox{}+\mbb{E}^{x_0} \left( \int_0^{t \wedge \tau}\!\! \smash[b]{\lint_{|y| < 1}} g_{R/4}(X_s^{(2)}) \left(\varrho(X_s^{(1)}+y) -\varrho(X_s^{(1)})-\nabla \varrho(X_{s}^{(1)}) \cdot y\right)  \nu(X_{s-},dy) \, ds \right). \notag
    \end{align}
	Using the submultiplicativity of $g \geq 0$ and \eqref{erg-eq24}, we find that the first term on the right-hand side is bounded above by
    \begin{align*}
        \mbb{E}^{x_0}& \left( \int_0^t \!\! \int_{|y| \geq 1} \varrho(X_s^{(1)}) \big[g_{R/16}(X_s^{(2)}) g(y) + g(X_s^{(2)}) g_{R/16}(y)\big] \, \nu(X_{s-},dy) \, ds \right) \\
        &\leq M_{R/16}(\mbb{R}^d) \int_0^t\mbb{E}^{x_0}(\varrho(X_s^{(1)}) g(X_s^{(2)})) \, ds + M(\mbb{R}^d) \int_0^t\mbb{E}^{x_0}(\varrho(X_s^{(1)}) g_{R/16}(X_s^{(2)}) \, ds
	\end{align*}
    with $M(\mbb{R}^d)$ and $M_{R/16}(\mbb{R}^d)$ from Definition~\ref{erg-11}. On the other hand, a similar calculation as in the proof of \eqref{erg-eq23} shows that the second term on the right-hand side of \eqref{erg-eq27} is less or equal than
    \begin{align*}
		C\mbb{E}^{x_0} \left( \int_0^t g_{R/4}(X_s^{(2)}) \varrho(X_s^{(1)}) \psi(X_s^{(1)}) \, ds \right)
	\end{align*}
    where $C$ is a suitable constant and $\psi \in C_b^2(\mbb{R}^d)$ such that $\spt \psi \cap \ball{0}{1/16}=\emptyset$. If we combine both estimates, apply Lemma~\ref{erg-15} and use that $\lim_{R\to\infty}M_{R/16}(\mbb{R}^d) = 0$, we get \eqref{erg-eq25}.

	In the general case, i.\,e.\ if $b(x) \neq 0$ or $Q(x) \neq 0$, we replace $X_t^{(1)}$ by
    \begin{equation*}
		X_t^{(1)} := \int_0^t b(X_{s}) \,ds + X_t^C +  \int_0^t \!\! \int_{0<|y| < 1} y  \, (N(dy,ds)-\mu(dy,ds))
	\end{equation*}
    where $(X_t^C)_{t \geq 0}$ denotes the continuous martingale part, cf.\ \cite[Theorem II.2.34]{js}; this gives additional terms when applying It\^o's formula, but the reasoning works exactly as in the pure-jump case.
\end{proof}

\begin{lemma} \label{erg-15}
    Let $(X_t)_{t \geq 0}$, $K$, $g$ and $x_0 \in \mbb{R}^d$ be as in Theorem~\ref{erg-13}. For any $T>0$ and all functions $g,\theta \in C_b^2(\mbb{R}^d)$ such that $\spt \theta \cap \ball{0}{\eps}=0$ for some sufficiently small $\eps>0$, there exists a constant $C>0$ such that
    \begin{align*}
		\mbb{E}^{x_0}\left(\exp \left[\sqrt{|X_{t \wedge \tau_K}^{(1)}|^2+1}-1 \right] g(X_{t \wedge \tau_K}^{(2)})\right) &\leq C \\
		\mbb{E}^{x_0}\left(\exp \left[\sqrt{|X_{t \wedge \tau_K}^{(1)}|^2+1}-1 \right] g(X_{t \wedge \tau_K}^{(2)}) \theta(X_t^{(1)})\right) &\leq Ct \\
		\mbb{E}^{x_0}\left(\exp \left[\sqrt{|X_{t \wedge \tau_K}^{(1)}|^2+1}-1 \right] g(X_{t \wedge \tau_K}^{(2)}) \theta(X_t^{(2)})\right) &\leq Ct
	\end{align*}
    for all $t \leq T$; here $\tau_K$ denotes the exit time from the set $K$ and $X_t-x_0 = X_t^{(1)}+X_t^{(2)}$ the decomposition from the proof of Theorem~\ref{erg-13}.
\end{lemma}
\begin{proof}
	We know from the proof of \cite[Theorem 4.1]{kuehn} that under the assumptions of Theorem~\ref{erg-13}
    \begin{equation*}
	   \sup_{t \leq T}\mbb{E}^{x_0}\left(\exp \left[\sqrt{|X_{t \wedge \tau_K}^{(1)}|^2+1}-1 \right] g(X_{t \wedge \tau_k}^{(2)}) \right)<\infty,
	\end{equation*}
    and this proves the first assertion. The other two estimates now follow from a straightforward application of It\^o's formula; mind that the initial term
    \begin{equation*}
		\exp \left[\sqrt{|X_{t \wedge \tau_K}^{(1)}|^2+1}-1 \right] g(X_{t \wedge \tau_K}^{(2)}) \theta(X_t^{(i)}) \Big|_{t=0} = 0
	\end{equation*}
	vanishes for $i \in \{1,2\}$ since $\theta(X_0^{(i)})=0$.
\end{proof}

\begin{remark}
    (i).
    The proof of Theorem~\ref{erg-13} simplifies substantially if the submultiplicative function $g \in C^2(\mbb{R}^d)$ satisfies the inequality \begin{equation}\label{erg-eq29}
    		|\nabla^2 g(x)| \leq C |g(x)|, \qquad x \in \mbb{R}^d,
    \end{equation}
    for some absolute constant $C>0$. In this case, we can apply It\^o's formula directly to the mapping $x \mapsto g(x-x_0) \chi_R(x-x_0)$ to prove \eqref{erg-eq21}; there is no need to use the decomposition $X_t = x+X_t^{(1)}+X_t^{(2)}$ and estimate \eqref{erg-eq22}. Although there are many examples of submultiplicative functions satisfying \eqref{erg-eq29}, it does not hold true for all (twice differentiable) submultiplicative functions.

    \medskip\noindent (ii).
    In Theorem~\ref{erg-13} submultiplicativity of the dominating function $g$ is required. This assumption can be weakened; it suffices to assume that there exist a subadditive function $a: \mbb{R}^d \to \mbb{R}$ and a submultiplicative function $m: \mbb{R}^d \to (0,\infty)$ such that $g(x) = m(x) \cdot a(x)$ for all $x \in \mbb{R}^d$, $a,m \in C^2(\mbb{R}^d)$ and
    \begin{equation*}
    		\lim_{R \to \infty} \inf_{|x| \geq R} |a(x)|>0.
    \end{equation*}
    The proof of Theorem~\ref{erg-13} under this relaxed assumption is similar, but more technical.
\end{remark}

Using exactly the same reasoning as in the proof of Theorem~\ref{erg-13}, we obtain a similar extension of Theorem~\ref{erg-10} to unbounded functions.
\begin{theorem} \label{erg-17}
    Let $(X_t)_{t \geq 0}$ be a rich L\'evy-type process with characteristics $(b,0,\nu)$ and symbol $q$, and let $f:\mbb{R}^d \to \mbb{R}$ be a continuous function satisfying the growth condition \eqref{G}. Suppose that either \eqref{B1} or \eqref{B2} holds for some $x_0 \in K$ and that $q(x_0,\cdot)$ satisfies the sector condition. Then the limit
    \begin{equation*}
        \lim_{t \to 0} \frac{1}{t} \left(\mbb{E}^{x_0} f(X_{t \wedge \tau_K})-f(x_0)\right)
    \end{equation*}
    exists and equals $Lf(x_0)$ as defined in Theorem~\ref{erg-10}. If $(X_t)_{t \geq 0}$ has bounded coefficients, then $K=\mbb{R}^d$ is admissible.
\end{theorem}

We close this section with an application of Corollary~\ref{erg-7}, which has been announced (without proof) in the recent publication \cite[remark following Theorem 5.2]{kuehn} on moments of L\'evy-type processes.
\begin{proposition} \label{app-1}
    Let $(X_t)_{t \geq 0}$ be a rich L\'evy-type process with symbol $q$ and characteristics $(b,Q,\nu)$. If there exist $x \in \mbb{R}^d$, $R \geq 0$ and $\alpha>0$ such that
    \begin{equation*}
	   \liminf_{t \to 0} \frac{1}{t} \mbb{E}^x\left(|X_t-x|^{\alpha} \I_{\{|X_t-x| > R\}}\right) < \infty,
	\end{equation*}
	then
    \begin{equation*}
		\int_{|y| > R} |y|^{\alpha} \, \nu(x,dy)
        \leq R^{\alpha} \nu(x,\{y \in \mbb{R}^d; |y| > R\}) \I_{R>0} + \liminf_{t \to 0} \frac{1}{t} \mbb{E}^x(|X_t-x|^{\alpha} \I_{\{|X_t-x| > R\}});
	\end{equation*}
	in particular $\int_{|y| > R} |y|^{\alpha} \, \nu(x,dy)<\infty$.
\end{proposition}

For $R=0$ Proposition~\ref{app-1} shows
\begin{equation*}
    C:=\liminf_{t \to 0} \frac{1}{t} \mbb{E}^x(|X_t-x|^{\alpha}) < \infty
    \implies \lint  |y|^{\alpha} \, \nu(x,dy) \leq C <\infty.
\end{equation*}

\begin{proof}[Proof of Proposition~\ref{app-1}]
	Since the identity
    \begin{equation}\tag{$\star$}\label{app-st3}
		\int |y|^{\alpha} \, \mu(dy)
        = \alpha \int_{(0,\infty)} \mu(|y| \geq r) r^{\alpha-1} \, dr
	\end{equation}
	holds for any $\alpha>0$ and any $\sigma$-finite measure $\mu$, we have
    \begin{equation*}
		\int_{|y| > R} |y|^{\alpha} \, \nu(x,dy)
		= \alpha \int_{(0,\infty)} \nu(x, \{y \in \mbb{R}^d; |y| > R, \: |y| \geq r\}) \, r^{\alpha-1} \, dr.
	\end{equation*}
	If $R=0$ then it follows from \eqref{erg-eq5} and Fatou's lemma that
    \begin{equation*}
		\int_{|y|>0} |y|^{\alpha} \, \nu(x,dy)
		\leq \alpha \liminf_{t \to 0} \frac{1}{t} \int_{(0,\infty)} \mbb{P}^x(|X_t-x| \geq r) \, r^{\alpha-1} \, dr
		\stackrel{\eqref{app-st3}}{=} \liminf_{t \to 0} \frac{1}{t} \mbb{E}^x(|X_t-x|^{\alpha}).
	\end{equation*}
    Here we use that the $\sigma$-finiteness of $\nu(x,dy)$ implies $\nu(x,\partial \ball{0}{r})=0$ for Lebesgue-almost all $r>0$. If $R>0$, then we split the integral
    \begin{equation*}
		\int_{|y| > R} |y|^{\alpha} \, \nu(x,dy)
		\leq R^{\alpha} \nu(x,\{y \in \mbb{R}^d; |y| > R\}) + \alpha \int_{(R,\infty)} \nu(x, \{y \in \mbb{R}^d; |y| \geq r\})\, r^{\alpha-1} \, dr,
	\end{equation*}
	and use again \eqref{erg-eq5} and Fatou's lemma to estimate the second term.
\end{proof}

\section{Uniform limits} \label{s-app}

In the previous section we have seen that the pointwise limit $\lim_{t \to 0} t^{-1} (\mbb{E}^{x_0} f(X_t)-f(x_0))$ exists for some fixed $x_0 \in \mbb{R}^d$ if $f \in C_0(\mbb{R}^d)$ satisfies a H\"{o}lder condition at $x_0$. Now we turn to the question under which assumptions on the regularity of $f$ the limit
\begin{equation}\label{app-eq0}
	\lim_{t \to 0} \frac 1t \left(\mbb{E}^{\bullet} f(X_t)-f(\bullet)\right)
\end{equation}
exists uniformly in $C_0(\mbb{R}^d)$, i.\,e.\ under which assumptions $f$ is contained in the domain $\mc{D}(A)$ of the generator of $(X_t)_{t \geq 0}$. It is well known that the limit exists (uniformly) for any function $f \in C_0^2(\mbb{R}^d)$ and any L\'evy-type process $(X_t)_{t \geq 0}$ with bounded coefficients, cf.\ \cite[Theorem~2.37]{ltp}. However, the results from the previous section suggest that the uniform limit may also exist for functions whose regularity varies from point to point, e.\,g.\ functions which satisfy
\begin{equation*}
	|f(x+y)-f(x)| \leq C |y|^{\alpha(x)} \fa x,y \in \mbb{R}^d,\; |y| \leq 1
\end{equation*}
for some absolute constant $C>0$ and a suitable mapping $\alpha: \mbb{R}^d \to [0,2]$. In this section, we will show that this is indeed true; more precisely we will establish that certain H\"{o}lder spaces of variable order are contained in the domain of the generator, cf.\ Corollary~\ref{app-11} and Corollary~\ref{app-13}. The idea is to use the fact that for a L\'evy-type process $(X_t)_{t \geq 0}$ the limit \eqref{app-eq0} exists uniformly if, and only if, the pointwise limit exists for each $x \in \mbb{R}^d$ and the limit defines a function in $C_0(\mbb{R}^d)$, cf.\ \cite[Theorem 7.22]{bm2}. At the end of this section we will present some examples, including stable-like and relativistic stable-like processes.

\medskip
Our first main result, Theorem~\ref{app-3}, is about the particular case that $(X_t)_{t \geq 0}$ is a L\'evy process.

\begin{proof}[Proof of Theorem~\ref{app-3}]
    \phantomsection\label{proof-app-3}\ref{app-3-i} is well known, see e.\,g.\ \cite[Theorem 31.5]{sato} or \cite[Theorem 2.37]{ltp}. The proofs of \ref{app-3-ii} and \ref{app-3-iii} are very similar, and therefore we only prove \ref{app-3-ii}. Pick a cut-off function $\chi \in C_c^{\infty}(\mbb{R}^d)$ such that $\chi \geq 0$, $\spt \chi \subseteq \ball{0}{1}$ and $\int_{\mbb{R}^d} \chi(x) \, dx = 1$. If we set $\chi_{\eps}(x) := \eps^{-1} \chi (\eps^{-1} x)$, then the convolution $f_n := \chi_{1/n} \ast f$ is in $C_0^2(\mbb{R}^d)$, hence in $\mc{D}(A)$, and $\lim_{n\to\infty}\|f_n-f\|_{\infty} = 0$. As
    \begin{align*}
		|(f_n-f)(x+y)-(f_n-f)(x)|
		&\leq \left| \int \chi(z) (f(x+y+n^{-1} z)-f(x+n^{-1} z)) \, dz \right| \\
		&\quad + |f(x+y)-f(x)| \\
		&\leq 2\|f\|_{\alpha} |y|^{\alpha}
	\end{align*}
	for all $|y| \leq 1$ and
    \begin{equation*}
		|(f_n-f)(x+y)-(f_n-f)(x)| \leq 2 \sup_{|r-s| \leq n^{-1}} |f(r)-f(s)| \xrightarrow[]{n \to \infty} 0,
	\end{equation*}
	we find
    \begin{equation*}
		\sup_{x \in \mbb{R}^d} \sup_{0<|y| \leq 1} \frac{|(f_n-f)(x+y)-(f_n-f)(x)|}{|y|^{\alpha}} \xrightarrow[]{n \to \infty} 0
	\end{equation*}
	which implies that
    \begin{equation*}
		Af_n(x) = \lint \left(f_n(x+y)-f_n(x)\right)  \nu(dy) \xrightarrow[]{n \to \infty} \lint \left(f(x+y)-f(x)\right)  \nu(dy)
	\end{equation*}
	uniformly in $x \in \mbb{R}^d$. Since the generator $(A,\mc{D}(A))$ is a closed operator, this finishes the proof.
\end{proof}

Next we extend Theorem~\ref{app-3} to L\'evy-type processes, cf.\ Theorem~\ref{app-5}.

\begin{proof}[Proof of Theorem~\ref{app-5}]
	\phantomsection\label{proof-app-5} It follows from Theorem~\ref{erg-9} that the pointwise limit
	\begin{equation*}
	Lf(x) = \lim_{t \to 0} \frac 1t\left(\mbb{E}^x f(X_t)-f(x)\right)
	\end{equation*}
	exists for all $x \in \mbb{R}^d$ and is given by
	\begin{itemize}
		\item
		$Lf(x) = \lint \left(f(x+y)-f(x)\right)  \nu(x,dy)$ for any $x \in \{0<\alpha \leq 1\}$;
		\item
		$Lf(x)= b(x) \cdot \nabla f(x) + \lint \left(f(x+y)-f(x)-\nabla f(x) \cdot y \I_{(0,1)}(|y|)\right)  \nu(x,dy)$ for any \\ $x \in \{1<\alpha<2\}$;
		\item
		$\begin{aligned}[t]\textstyle
		Lf(x) &= b(x) \cdot \nabla f(x) + \tfrac{1}{2} \tr\left(Q(x) \nabla^2 f(x)\right) \\ &\quad \textstyle + \lint \left(f(x+y)-f(x)-\nabla f(x) \cdot y \I_{(0,1)}(|y|)\right)  \nu(x,dy)\end{aligned}$ \\ for any $x \in \{\alpha=2\}$.
	\end{itemize}
	As $Q(x)=0$ for all $x \in \{0<\alpha<2\}$ and $\lint_{|y|<1} y \, \nu(x,dy) = b(x)$ for all $x \in \{0<\alpha \leq 1\}$ (see Lemma~\ref{app-7} in the appendix), we can write $Lf$ in a closed form as
	\begin{equation*}
	Lf(x) = b(x) \cdot g(x) +\frac{1}{2} \tr \left(Q(x) h(x)\right) + \lint \left(f(x+y)-f(x)-g(x) \cdot y \I_{(0,1)}(|y|)\right)  \nu(x,dy).
	\end{equation*}
	In order to prove that $f$ is contained in the domain of the generator $A$ and $Af = Lf$, it suffices to show that $Lf \in C_0(\mbb{R}^d)$, see e.\,g.\ \cite[Theorem 7.22]{bm2}. The triangle inequality, Taylor's formula and conditions \eqref{C1}-\eqref{C3} imply that there exists a constant $C>0$ such that
	\begin{equation}\label{app-eq9}
	|f(x+y)-f(x)-g(x) \cdot y| \leq C |y|^{\alpha(x)} \fa x,y \in \mbb{R}^d,\; |y| \leq 1.
	\end{equation}
	Fix a cut-off function $\chi \in C_c^{\infty}(\mbb{R}^d)$ such that $\chi \geq 0$, $\spt \chi \subseteq \ball{0}{1}$ and $\int_{\mbb{R}^d} \chi(x) \, dx = 1$. If we set $\chi_{\eps}(x) := \eps^{-1} \chi (\eps^{-1} x)$, then the convolutions $f_n := \chi_{1/n} \ast f$, $g_n := \chi_{1/n} \ast g$ and $h_n := \chi_{1/n} \ast h$ are $C_0^2(\mbb{R}^d)$-functions and
	\begin{equation*}
	\|f_n-f\|_{\infty}+\|g_n-g\|_{\infty}+\|h_n-h\|_{\infty} \xrightarrow[]{n \to \infty} 0.
	\end{equation*}
	We are going to show that
	\begin{equation*}
	\Delta_n(x,y) := (f_n-f)(x+y)-(f_n-f)(x)-(g_n-g)(x) \cdot y
	\end{equation*}
	satisfies an estimate similar to \eqref{app-eq9}. By the very definition of the convolution, we have
	\begin{align*}
	\Delta_n(x,y)
	&=  \int (f(x+y+z)-f(x+y)) \chi_{1/n}(z) \, dz - \int (f(x+z)-f(x)) \chi_{1/n}(z) \, dz\\
	&\quad - \int (g(x+z)-g(x)) \cdot y \chi_{1/n}(z) \, dz.
	\end{align*}
	Since $\spt \chi_{1/n} \subseteq \cball{0}{1/n}$ and $0 \leq \chi \leq 1$,
	\begin{equation*}
	|\Delta_n(x,y)| \leq 2\sup_{|r-s| \leq n^{-1}} |f(r)-f(s)|+ \sup_{|r-s| \leq n^{-1}} |g(r)-g(s)|.
	\end{equation*}
	On the other hand, we have by \eqref{app-eq9}
	\begin{align*}
	|\Delta_n(x,y)|
	\leq 2C \sup_{|x-z|\leq {n^{-1}}} |y|^{\alpha(z)} \int_{\ball{0}{1}} \chi(z) \, dz
	= 2C |y|^{\alpha(x)} \sup_{ |x-z|\leq {n^{-1}}} |y|^{\alpha(z)-\alpha(x)}.
	\end{align*}
	As $\alpha$ is uniformly continuous, we can choose $N \in \mbb{N}$ sufficiently large such that \begin{equation*}
	|\alpha(x)-\alpha(z)| \leq \eps/2 \fa x \in \mbb{R}, z \in \ball{x}{N^{-1}}.
	\end{equation*}
	Combining both estimates, we find
	\begin{align*}
	&\frac{|\Delta_n(x,y)|}{|y|^{\alpha(x)-\eps}}\\
	&\leq \min \left\{ \frac{2\sup_{|r-s| \leq n^{-1}} |f(r)-f(s)|+ \sup_{|r-s| \leq n^{-1}} |g(r)-g(s)|}{|y|^{\alpha(x)-\eps}}, \: 2C \sup_{|x-z|\leq {n^{-1}}} |y|^{\eps+(\alpha(z)-\alpha(x))} \right\} \\
	&\leq \min \left\{  \frac{2\sup_{|r-s| \leq n^{-1}} |f(r)-f(s)|+ \sup_{|r-s| \leq n^{-1}} |g(r)-g(s)|}{|y|^{2}}, \: 2C |y|^{\eps/2} \right\}
	\end{align*}
	for all $x \in \mbb{R}^d$, $0<|y| \leq 1$ and $n \geq N$. As $f \in C_0(\mbb{R}^d)$ and $g \in C_0(\mbb{R}^d)$ are uniformly continuous, this proves 
	\begin{align*}
	\lim_{n \to \infty} \sup_{x \in \mbb{R}^d} \sup_{0<|y| \leq 1} \frac{|\Delta_n(x,y)|}{|y|^{\alpha(x)-\eps}} =0.
	\end{align*}
	In particular, there exist constants $C_n>0$ such that $C_n \to 0$ as $n \to \infty$ and
	\begin{align*}
	|(f_n-f)(x+y)-(f_n-f)(x)-(g_n-g)(x) \cdot y| &\leq C_n |y|^{\alpha(x)-\eps}
	\end{align*}
	for all $x,y \in \mbb{R}^d$, $|y| \leq 1$. If we set
	\begin{align*}
	Lf_n(x) 
	&:= b(x) g_n(x) + \frac{1}{2} \tr\left(Q(x) h_n(x)\right) \\
	&\quad + \lint \left(f_n(x+y)-f_n(x)-g_n(x) \cdot y \I_{(0,1)}(|y|) \right)  \nu(x,dy),
	\end{align*}
	then
	\begin{align*}
	|L f_n(x)- Lf(x)|
	&\leq \|b\|_{\infty} \|g_n-g\|_{\infty} + \|Q\|_{\infty} \|h_n-h\|_{\infty}+ C_n \lint_{|y| \leq 1} |y|^{\alpha(x)-\eps} \, \nu(x,dy) \\
	&\quad +2 \|f_n-f\|_{\infty} \sup_{x \in \mbb{R}^d} \int_{|y|>1} \nu(x,dy).
	\end{align*}
	This expression converges to zero uniformly in $x$ since $(X_t)_{t \geq 0}$ has bounded coefficients. As $Lf_n \in C_0(\mbb{R}^d)$ for large $n \in \mbb{N}$, see Lemma~\ref{app-9} below, we conclude that $Lf \in C_0(\mbb{R}^d)$.
\end{proof}	

For the proof of Theorem~\ref{app-5} we need the following auxiliary statement.
\begin{lemma} \label{app-9}
	$Lf_n$ defined in the proof of Theorem~\ref{app-5} is a $C_0(\mbb{R}^d)$-function for sufficiently large $n \in \mbb{N}$.
\end{lemma}

\begin{proof}
	The mapping $x \mapsto Q(x)$ is, by assumption, continuous and bounded. As $h_n \in C_0^2(\mbb{R}^d)$, this implies that $\tr(Q(\bullet)^\top h_n(\bullet)) \in C_0(\mbb{R}^d)$. Consequently, it is enough to show that
	\begin{equation*}
	\tilde{L}f_n(x) :=  b(x) \cdot g_n(x)+ \lint \left(f_n(x+y)-f_n(x)-g_n(x) \cdot y \I_{(0,1)}(|y|) \right)  \nu(x,dy) \in C_0(\mbb{R}^d).
	\end{equation*}
	Since $C_c^{\infty}(\mbb{R}^d) \subseteq \mc{D}(A)$ and $(X_t)_{t \geq 0}$ has bounded coefficients, we have $C_0^2(\mbb{R}^d) \subseteq \mc{D}(A)$, and therefore
	\begin{align*}
	Af_n(x)
	&= b(x) \cdot \nabla f_n(x) + \frac{1}{2} \tr\left(Q(x) \nabla^2 f_n(x)\right) \\
	&\quad + \lint \left(f_n(x+y)-f_n(x)-\nabla f_n(x) \cdot y \I_{(0,1)}(|y|) \right)  \nu(x,dy) 
	\end{align*}
	is in $C_0(\mbb{R}^d)$. Using again the fact that $Q \in C_b(\mbb{R}^d)$ and $\nabla^2 f_n \in C_0(\mbb{R}^d)$, we get
	\begin{equation*}
	\tilde{A}f_n(x) := b(x) \cdot \nabla f_n(x) + \lint \left(f_n(x+y)-f_n(x)-\nabla f_n(x) \cdot y \I_{(0,1)}(|y|) \right)  \nu(x,dy) \in C_0(\mbb{R}^d).
	\end{equation*}
	Let $x \in \mbb{R}^d$. We distinguish between two cases.
	\begin{description}
		\item[\normalfont $0<\alpha(x) \leq 1+\eps/2$:]
		Using our assumption $\beta_{\infty}^{x} + \eps \leq \alpha(x)$, we find $\beta_{\infty}^{x}<1$ which implies, by Lemma~\ref{app-7}, $b(x) - \int_{|y|<1} y \, \nu(x,dy)=0$. Thus, $\tilde{A}f_n(x) = \tilde{L}f_n(x)$.
		\item[\normalfont $1+\eps/2 <\alpha(x)$:]
		Since $\alpha$ is uniformly continuous, we can choose $n \in \mathbb{N}$ (not depending on $x$) so large that $|\alpha(x)-\alpha(z)| \leq \eps/4$ for all $z \in \cball{x}{n^{-1}}$. Then $\alpha(z)>1+\eps/4$ for all $z \in \ball{x}{n^{-1}}$ and, therefore, $f|_{\ball{x}{n^{-1}}}$ is differentiable. As $\spt \chi_{1/n} \subseteq \cball{0}{1/n}$, this implies $\nabla f_n(x) = g_n(x)$. Hence, $\tilde{L}f_n(x) = \tilde{A}f_n(x)$.
	\end{description}
	Consequently, we have $\tilde{L}f_n = \tilde{A}f_n \in C_0(\mbb{R}^d)$ for $n \in \mbb{N}$ sufficiently large.
\end{proof}

\begin{corollary} \label{app-11}
    Let $(X_t)_{t \geq 0}$ be a rich L\'evy-type process with symbol $q$, $q(x,0)=0$ and characteristics $(b,0,\nu)$. Suppose that $q$ has bounded coefficients and $b(x)  =\lint_{|y|<1} y \, \nu(x,dy)$ for all $x \in \mbb{R}^d$. Let $\eps>0$ and $\alpha: \mbb{R}^d \to [\eps,1]$ be uniformly continuous such that
    \begin{equation*}
		\sup_{x \in \mbb{R}^d} \lint_{|y| \leq 1} |y|^{\alpha(x)-\eps} \, \nu(x,dy)<\infty.
	\end{equation*}
    If either the sector condition \eqref{app-eq11} holds or $\beta_{\infty}^x \leq \alpha(x)-\eps$ for all $x \in \mbb{R}^d$, then the H\"{o}lder space of variable order
    \begin{equation*}
        \mc{C}_0^{\alpha(\cdot)}
        := \left\{f \in C_0(\mbb{R}^d);\: \sup_{x \in \mbb{R}^d} \sup_{0<|y| \leq 1} \frac{|f(x+y)-f(x)|}{|y|^{\alpha(x)}}<\infty \right\}
	\end{equation*}
	is contained in the domain of the generator $A$ and
    \begin{equation*}
		Af(x) = \lint \left(f(x+y)-f(x)\right)  \nu(x,dy)
        \fa x \in \mbb{R}^d,\; f \in \mc{C}_0^{\alpha(\cdot)}.
	\end{equation*}
\end{corollary}
\begin{proof}
    Under the assumptions of Corollary~\ref{app-11}, we know from the remark following Theorem~\ref{app-5} that $\beta_{\infty}^{x} \leq \alpha(x)-\eps$ for all $x \in \mbb{R}^d$. Moreover, $\alpha(x) \in [0,1]$ for all $x \in \mbb{R}^d$ and, by assumption, condition \eqref{C1} is satisfied for all $x \in \mbb{R}^d$. Consequently, the assumptions of Theorem~\ref{app-5} are satisfied, and so Theorem~\ref{app-5} proves the assertion.
\end{proof}

Let us mention that among the first to consider H\"{o}lder spaces of variable order were Ross \& Samko \cite{ross} who study fractional integrals of variable order. In \cite{almeida} H\"{o}lder spaces of variable order are shown to be particular cases of Besov spaces with variable smoothness and integrability; see Andersson \cite{andersson} for further characterizations.

\begin{corollary}\label{app-13}
    Let $(X_t)_{t \geq 0}$ be a rich L\'evy-type process with bounded coefficients and with symbol $q$ and characteristics $(b,0,\nu)$. Let $\eps>0$ be a constant and $\alpha: \mbb{R}^d \to (\eps,2]$ be a uniformly continuous mapping. Suppose that either the sector condition \eqref{app-eq11} is satisfied or $\alpha(x)-\eps \geq \beta_{\infty}^x$ for all $x \in \mbb{R}^d$. If
    \begin{equation*}
		\sup_{x \in \mbb{R}^d} \lint_{|y| \leq 1} |y|^{\alpha(x)-\eps} \, \nu(x,dy)<\infty,
	\end{equation*}
	then the space 
    \begin{equation*}
		\mc{C}_0^{1,(\alpha(\cdot)-1)^+}	
        :=  \left\{f \in C^1_0(\mbb{R}^d); \forall j=1,\ldots,d: \partial_j f \in \mc{C}_0^{\max\{\alpha(\cdot)-1,0\}}\right\}
	\end{equation*}
	is contained in the domain of the generator $A$, and for all $x \in \mbb{R}^d$ and $f \in \mc{C}_0^{1,(\alpha(\cdot)-1)^+}$
    \begin{equation*}
        Af(x) = b(x) \cdot \nabla f(x) + \lint \left(f(x+y)-f(x)-\nabla f(x) \cdot y \I_{(0,1)}(|y|)\right)  \nu(x,dy) .
	\end{equation*}
\end{corollary}
\begin{proof}
    As $\mc{C}_0^{1,(\alpha(\cdot)-1)^+} \subseteq C_0^1(\mbb{R}^d)$, we may assume without loss of generality that $\alpha(x) \geq 1$ for all $x \in \mbb{R}^d$; otherwise we could replace $\alpha$ by $\max\{\alpha,1\}$. As in the proof of Corollary~\ref{app-11}, we find $\beta_{\infty}^{x} +\eps \leq \alpha(x)$ for all $x \in \mbb{R}^d$.  It remains to check that $f \in \mc{C}_0^{1,\alpha(\cdot)-1}$ satisfies the assumptions of Theorem~\ref{app-5}. If $x \in \mbb{R}^d$ is such that $\alpha(x)=1$ it is obvious from the mean value theorem that \eqref{C1} is satisfied. Now let $x \in \{1<\alpha<2\}$. Applying the mean value theorem to the auxiliary function $h(y) := f(x+y)-f(x)-\nabla f(x) \cdot y$ shows
    \begin{align*}
        |f(x+y)-f(x)-\nabla f(x) \cdot y|
        \leq |y| \sup_{\zeta \in \ball{0}{|y|}} |\nabla f(x+\zeta)-\nabla f(x)| \leq C |y| \cdot |y|^{\alpha(x)-1}
        = C |y|^{\alpha(x)}
	\end{align*}
    for all $x,y \in \mbb{R}^d$, $|y| \leq 1$ and some absolute constant $C>0$; here we use that $\partial_j f \in \mc{C}_0^{\alpha(\cdot)-1}$ for all $j \in \{1,\ldots,d\}$. This shows that condition \eqref{C2} holds true.
\end{proof}

We close this section with some examples. Recall the definition of the H\"{o}lder spaces of variable order $\mc{C}_0^{\alpha(\cdot)}$ and $\mc{C}_0^{1,(\alpha(\cdot)-1)^+}$ introduced in Corollary~\ref{app-11} and Corollary~\ref{app-13}, respectively.

\begin{example}[Stable-like dominated process] \label{app-15}
    Let $(X_t)_{t \geq 0}$ be a rich L\'evy-type process with symbol $q$ and characteristics $(b,0,\nu)$. Denote by $(A,\mc{D}(A))$ the generator of $(X_t)_{t \geq 0}$. Suppose that $(X_t)_{t \geq 0}$ has bounded coefficients and that there exist a constant $c>0$ and a mapping $\gamma:\mbb{R}^d \to (0,2)$ such that $\inf_{x \in \mbb{R}^d} \gamma(x)>0$ and
    \begin{equation*}
        \nu(x,A \cap \ball{0}{1}) \leq c \int_{A \cap \ball{0}{1}} \frac{dy}{|y|^{d+\gamma(x)}}
        \fa A \in \mc{B}(\mbb{R}^d\setminus \{0\}),\; x \in \mbb{R}^d.
	\end{equation*}
    Let $\alpha: \mbb{R}^d \to (0,2)$ be a uniformly continuous mapping such that $\inf_{x \in \mbb{R}^d} (\alpha(x)-\gamma(x))>0$, and suppose that either the sector condition \eqref{app-eq11} is satisfied or $\inf_{x \in \mbb{R}^d} (\alpha(x)-\beta_{\infty}^x)>0$.
    \begin{enumerate}
        \item
            If $\alpha(\mbb{R}^d) \subseteq [0,1]$ and $b(x) = \lint_{|y|<1} y \, \nu(x,dy)$ for all $x \in \mbb{R}^d$, then $\mc{C}_0^{\alpha(\cdot)} \subseteq \mc{D}(A)$ and
            \begin{equation*}
                Af(x) = \lint  \left(f(x+y)-f(x)\right)  \nu(x,dy),
                \qquad x \in \mbb{R}^d,\; f \in \mc{C}_0^{\alpha(\cdot)}.
		    \end{equation*}
		\item
            $\mc{C}_0^{1,(\alpha(\cdot)-1)^+} \subseteq \mc{D}(A)$ and
            \begin{equation*}
                Af(x) = b(x) \cdot \nabla f(x)+ \lint  \left(f(x+y)-f(x)-\nabla f(x) \cdot y \I_{(0,1)}(|y|)\right)  \nu(x,dy)
		    \end{equation*}
		for all $f \in \mc{C}_0^{1,(\alpha(\cdot)-1)^+}$ and $x \in \mbb{R}^d$.
	\end{enumerate}
\end{example}

\begin{example} \label{app-17}
	Let $(X_t)_{t \geq 0}$ be a rich L\'evy-type process with one of the following symbols.
    \begin{itemize}
	\item
        stable-like: $q(x,\xi) = |\xi|^{\gamma(x)}$ where $\gamma: \mbb{R}^d \to (0,2)$ is a H\"{o}lder continuous mapping such that $\inf_{x \in \mbb{R}^d} \gamma(x)>0$.
	
    \item
        relativistic stable-like: $q(x,\xi) = (|\xi|^2+m(x)^2)^{\gamma(x)/2}-m(x)^{\gamma(x)}$ for H\"{o}lder continuous mappings $\gamma: \mbb{R}^d \to (0,2)$ and $m: \mbb{R}^d \to (0,\infty)$ such that
        \begin{equation*}
		    \inf_{x \in \mbb{R}^d} \gamma(x)>0
            \quad\text{and}\quad
            0 < \inf_{x \in \mbb{R}^d} m(x) \leq \sup_{x \in \mbb{R}^d} m(x) < \infty.
	    \end{equation*}

	\item
        TLP-like:\footnote{TLP is short for ``truncated L\'evy process''.}
        $q(x,\xi) = (|\xi|^2+m(x)^2)^{\gamma(x)/2} \cos\big[\gamma(x) \arctan \tfrac{|\xi|}{m(x)}\big]-m(x)^{\gamma(x)}$ for H\"{o}lder continuous mappings $\gamma: \mbb{R}^d \to (0,1)$ and $m: \mbb{R}^d \to (0,\infty)$ such that
        \begin{equation*}
            0<\inf_{x \in \mbb{R}^d} \gamma(x) \leq \sup_{x \in \mbb{R}^d} \gamma(x)<1
            \quad\text{and}\quad
            0 < \inf_{x \in \mbb{R}^d} m(x) \leq \sup_{x \in \mbb{R}^d} m(x) < \infty.
		\end{equation*}
		
    \item
        Lamperti stable-like: $q(x,\xi) = (|\xi|^2+m(x))_{\gamma(x)}-(m(x))_{\gamma(x)}$ -- $(z)_{\gamma} := \Gamma(z+\gamma)/\Gamma(z)$ denotes the Pochhammer symbol -- for H\"{o}lder continuous mappings $\gamma: \mbb{R}^d \to (0,1)$ and $m: \mbb{R}^d \to (0,\infty)$ such that
        \begin{equation*}
            0<\inf_{x \in \mbb{R}^d} \gamma(x) \leq \sup_{x \in \mbb{R}^d} \gamma(x)<1
            \quad \text{and} \quad
            0 < \inf_{x \in \mbb{R}^d} m(x) \leq \sup_{x \in \mbb{R}^d} m(x) < \infty.
		\end{equation*}
	\end{itemize}
	Let $\alpha: \mbb{R}^d \to [0,2]$ be a uniformly continuous mapping such that $\inf_{x \in \mbb{R}^d} (\alpha(x)-\gamma(x))>0$. Then:
    \begin{enumerate}
	   \item
            $\mc{C}_0^{1,(\alpha(\cdot)-1)^+} \subseteq \mc{D}(A)$ and
            $Af(x) = \lint \left(f(x+y)-f(x)-\nabla f(x) \cdot y \I_{(0,1)}(|y|)\right)  \nu(x,dy)$ for any $f \in \mc{C}_0^{1,(\alpha(\cdot)-1)^+}$ where $(0,0,\nu)$ denotes the characteristics of the symbol $q$.

		\item
            If $\alpha(\mbb{R}^d) \subseteq [0,1]$, then $\mc{C}_0^{\alpha(\cdot)} \subseteq \mc{D}(A)$ and $Af(x) = \lint \left(f(x+y)-f(x)\right)  \nu(x,dy)$ for any $f \in \mc{C}_0^{\alpha(\cdot)}$.
	\end{enumerate}
\end{example}

Example~\ref{app-17} is a direct consequence of Theorem~\ref{app-5} and Remark~\ref{app-6}(ii) since all symbols satisfy the sector condition \eqref{app-eq11} and growth condition \eqref{app-eq12} with $\delta = 0$. Note that the existence of (rich) L\'evy-type processes with the symbols mentioned in Example~\ref{app-17} has been established in \cite{diss} recently, see also \cite{matters}. Obviously, Example~\ref{app-17} applies, in particular, in the L\'evy case, i.\,e.\ if the maps $\gamma(\bullet)$ and $m(\bullet)$ are constants.

\medskip
We close this section with the following example.

\begin{example}[L\'evy-driven SDE]\label{app-19}
    Let $(L_t)_{t \geq 0}$ be a $k$-dimensional L\'evy process with L\'evy triplet $(0,0,\nu)$ and characteristic exponent $\psi$. Suppose that the L\'evy measure $\nu$ is symmetric and that there exists an $\alpha \in (0,2)$ such that $\lint_{|y| \leq 1} |y|^{\alpha} \, \nu(dy)<\infty$. For any bounded (globally) Lipschitz continuous function $\sigma: \mbb{R}^d \to \mbb{R}^{d \times k}$ the solution to the SDE
    \begin{equation*}
		dX_t = \sigma(X_{t-}) \, dL_t, \qquad X_0 = x,
	\end{equation*}
	is a rich L\'evy-type process with symbol $q(x,\xi) = \psi(\sigma(x)^\top \xi)$, $x,\xi \in \mbb{R}^d$. Moreover:,
    \begin{enumerate}
	   \item
            If $\alpha \in (0,1)$, then the H\"{o}lder space $\mc{C}_0^{\beta}$ is contained in the domain of the generator $A$ of $(X_t)_{t \geq 0}$ for any $\beta \in (\alpha,1]$ and
            \begin{equation*}
			    Af(x) = \lint  \left(f(x+\sigma(x)y)-f(x)\right)  \nu(dy),
                \qquad f \in \mc{C}_0^{\beta},\; x \in \mbb{R}^d.
		    \end{equation*}

        \item
            If $\alpha \in [1,2)$, then the H\"{o}lder space $\mc{C}_0^{1,\beta-1}$ is contained in the domain of the generator $A$ of $(X_t)_{t \geq 0}$ for any $\beta \in (\alpha,2]$ and
            \begin{equation*}
                Af(x) = \lint  \left(f(x+\sigma(x)y)-f(x)-\nabla f(x) \cdot \sigma(x) y \I_{(0,1)}(|y|)\right)  \nu(dy)
			\end{equation*}
			for any $f \in \mc{C}_0^{1,\beta-1}$ and $x \in \mbb{R}^d$.
\end{enumerate}
\end{example}

\begin{proof}
    It is well known that the solution $(X_t)_{t \geq 0}$ to the SDE is a rich L\'evy-type process with symbol $q(x,\xi) = \psi(\sigma(x)^\top \xi)$, cf.\  Schilling \& Schnurr \cite[Corollary 3.7]{schnurr} or K\"{u}hn \cite[Example 4.1]{sde}. Since the L\'evy measure $\nu$ is symmetric, both $\psi$ and $q$ are real-valued; in particular, $q$ satisfies the sector condition. Moreover, the characteristics of $q$ are given by $(0,0,\nu(x,dy))$ where
    \begin{equation*}
		\nu(x,B) := \int \I_B(\sigma(x) y) \, \nu(dy), \qquad x \in \mbb{R}^d,\; B \in \mc{B}(\mbb{R}^d \setminus \{0\});
	\end{equation*}
	therefore, the boundedness of $\sigma$ gives
    \begin{equation*}
		\sup_{x \in \mbb{R}^d} \lint_{|y| \leq 1} |y|^{\alpha} \, \nu(x,dy)
        \leq \|\sigma\|_{\infty}^{\alpha} \int_{|y| \leq 1} |y|^{\alpha} \, \nu(dy)
        < \infty.
	\end{equation*}
	Now the assertion follows from Corollary~\ref{app-11} and Corollary~\ref{app-13}.
\end{proof}

\section*{Acknowledgement}
    We thank the anonymous referees for their careful reading and helpful suggestions which helped to improve the presentation of our paper.

\appendix

\specialsection

\section{Appendix}

We frequently use Dynkin's formula which can be seen as a probabilistic counterpart of the fundamental theorem of calculus, see e.\,g.\ \cite[Lemma I.4.1.14]{jacob123} or \cite[Proposition 7.31]{bm2}.

\begin{lemma}[Dynkin's formula] \label{aux-0}
	Let $(X_t)_{t \geq 0}$ be a L\'evy-type process with infinitesimal generator $(A,\mc{D}(A))$, and let $x \in \mbb{R}^d$. If $\tau$ is a stopping time such that $\mbb{E}^x(\tau)<\infty$, then \begin{equation*}
		\mbb{E}^xf(X_{\tau})-f(x) = \mbb{E}^x \left( \int_{(0,\tau)} Af(X_s) \, ds \right)
	\end{equation*}
	for all $f \in \mc{D}(A)$ and $t \geq 0$.
\end{lemma}

Recall that a function $\psi: \mbb{R}^d \to \mbb{C}$ with $\psi(0)=0$ is \emph{continuous negative definite}, if it admits a L\'evy--Khintchine representation of the form \eqref{def-eq1}. A continuous negative definite function $\psi$ satisfies the \emph{sector condition} if there exists a constant $C>0$ such that
\begin{equation*}
	|\im \psi(\xi)| \leq C \re \psi(\xi) \fa \xi \in \mbb{R}^d.
\end{equation*}
The following lemma is used in the proof of Theorem~\ref{erg-10}.

\begin{lemma} \label{aux-1}
    Let $\psi$ be a continuous negative definite function with triplet $(b,0,\nu)$, and let $\alpha \in (0,2)$. The following statements are equivalent: \begin{enumerate}
		\item\label{aux-1-i} $\lint_{|y| \leq 1} |y|^{\alpha} \,\nu(dy)< \infty$;
		\item\label{aux-1-ii} $\int_1^{\infty} \sup_{|\xi| \leq r} \re \psi(\xi) \, \frac{dr}{r^{1+\alpha}} < \infty$;
		\item\label{aux-1-iii} $\int_{|\xi| \geq 1}  \re \psi(\xi) \, \frac{d\xi}{|\xi|^{d+\alpha}} < \infty$.
	\end{enumerate}
	If $\psi$ satisfies the sector condition, then we may replace $\re \psi$ by $|\psi|$.
\end{lemma}

\begin{proof}
    Obviously, it suffices to prove the first assertion. We prove $\ref{aux-1-i}\Rightarrow\ref{aux-1-ii} \Rightarrow \ref{aux-1-iii} \Rightarrow \ref{aux-1-i}$.

	$\ref{aux-1-i} \Rightarrow \ref{aux-1-ii}$: Since $1-\cos(y \cdot \xi) \leq \tfrac{1}{2} |y \xi|^2$ for all $y,\xi \in \mbb{R}^d$, we have 
	\begin{equation*}
		\psi^*(r) := \sup_{|\xi| \leq r} \re \psi(\xi) \leq 2 \lint  \min\left\{1,|y|^2 r^2\right\} \, \nu(dy)
	\end{equation*}
	implying
    \begin{equation*}
		\int_1^{\infty} \psi^*(r)\,\frac{dr}{r^{1+\alpha}}
		\leq 2 \int_1^{\infty} r^2 \lint_{|y| < r^{-1}} |y|^2 \, \nu(dy)\, \frac{dr}{r^{1+\alpha}}
            + 2 \int_1^{\infty} \int_{|y| \geq r^{-1}} \, \nu(dy) \,\frac{dr}{r^{1+\alpha}}
		=: 2I_1+2I_2.
	\end{equation*}
	An application of Tonelli's theorem shows
    \begin{align*}
		I_1
		= \int_1^{\infty} \lint_{|y| < r^{-1}} r^{1-\alpha} |y|^2 \, \nu(dy) \, dr
		&= \lint_{|y| \leq 1} \left( \int_{1 \leq r < |y|^{-1}} r^{1-\alpha} \, dr \right) |y|^2 \, \nu(dy) \\
		&= \frac{1}{2-\alpha} \lint_{|y| \leq 1} |y|^2 \left( |y|^{\alpha-2} -1 \right)  \nu(dy) <\infty
	\end{align*}
	and
    \begin{align*}
		I_2
		= \int_1^{\infty} \nu(\{y;\, |y|\geq r^{-1}\}) \,\frac{dr}{r^{1+\alpha}}
		&= \int_0^1 \nu(\{y;\, |y| \geq u\}) \,\frac{du}{u^{1-\alpha}}\\
		&= \frac 1\alpha \lint_{|y| \leq 1} |y|^{\alpha} \, \nu(dy) < \infty.
	\end{align*}
	In the last step we use the identity
    \begin{equation*}
		\int f(x) \, d\mu(x) = \int_0^{\infty} \mu(\{x;\, |f(x)| \geq r\}) \, dr
	\end{equation*}
    which holds for any $\sigma$-finite measure $\mu$ on $(\mbb{R}^d\setminus\{0\},\mc{B}(\mbb{R}^d \setminus \{0\}))$ and any non-negative measurable function $f$. This proves \ref{aux-1-ii}.

    The implication $\ref{aux-1-ii} \Rightarrow \ref{aux-1-iii}$ follows easily by introducing spherical coordinates and using the obvious estimate
    \begin{equation*}
		\re \psi(r \eta) \leq \psi^*(r)
        \fa r \geq 0,\; \eta \in \mbb{R}^d,\; |\eta|=1.
	\end{equation*}

	It remains to prove that \ref{aux-1-iii} implies \ref{aux-1-i}. To this end, we note that
    \begin{equation*}
		|y|^{\alpha} = c \int (1-\cos(y \cdot \xi)) \frac{1}{|\xi|^{d+\alpha}} \, d\xi
	\end{equation*}
	for the constant $c = \alpha 2^{\alpha-1} \pi^{-d/2} \Gamma\left(\frac{\alpha+d}{2}\right)\big/\Gamma\left(1-\frac\alpha 2\right)$.
	It is not difficult to see that this implies \begin{equation*}
		|y|^{\alpha} \leq c' \int_{|\xi| \geq 1} (1-\cos(y \cdot \xi)) \, \frac{1}{|\xi|^{d+\alpha}} \, d\xi \fa |y| \leq 1
	\end{equation*}
	for some constant $c'=c'(\alpha)>0$. Hence,
    \begin{equation}\label{aux-eq1}\begin{aligned}
        \lint_{|y| \leq 1} |y|^{\alpha} \, \nu(dy)
		&\leq c' \int_{|\xi| \geq 1} \left( \lint_{|y| \leq 1} \left(1-\cos(y \cdot \xi)\right)  \nu(dy) \right)\, \frac{d\xi}{|\xi|^{d+\alpha}}\\
		&\leq c' \int_{|\xi| \geq 1} \re \psi(\xi) \, \frac{d\xi}{|\xi|^{d+\alpha}} \, d\xi < \infty,
    \end{aligned}\end{equation}
	which completes the proof.
\end{proof}

For a continuous negative definite function $\psi$ the Blumenthal--Getoor index at $\infty$ can be defined by
\begin{equation*}
	\beta_{\infty} := \inf\left\{\gamma>0; \lim_{r \to \infty} \frac{1}{r^\gamma} \sup_{|\xi| \leq r} |\psi(\xi)| < \infty \right\},
\end{equation*}
cf.\ Schilling \cite{rs97} or Blumenthal \& Getoor \cite{blumen}. The following auxiliary statement is needed in the proof of Theorem~\ref{app-5}.

\begin{lemma} \label{app-7}
	Let $\psi$ be a continuous negative definite function,
    \begin{equation*}
		\psi(\xi) = ib \cdot \xi + \frac{1}{2} \xi \cdot Q \xi + \lint \left(1-e^{i\xi \cdot y} +i \xi \cdot y \I_{(0,1)}(|y|)\right)  \nu(dy),
        \qquad \xi \in \mbb{R}^d,
	\end{equation*}
	and denote by $\beta_{\infty} \in [0,2]$ the Blumenthal--Getoor index at $\infty$.
    \begin{enumerate}
		\item\label{app-7-i} If $\beta_{\infty} < 2$, then $Q=0$.
		\item\label{app-7-ii} If $\beta_{\infty}<1$, then $b = \int_{|y|<1} y \, \nu(dy)$.
	\end{enumerate}
\end{lemma}

\begin{proof}
    \ref{app-7-i}
    Since $|\xi|^{-2} |1-\cos(y \cdot \xi)| \leq \min\{2,|y|^2\}$ for all $|\xi| \geq 1$, an application of the dominated convergence theorem shows \begin{equation*}
    		\lim_{|\xi| \to \infty} \frac{1}{|\xi|^2} \lint  \left(1-\cos(y \cdot \xi)\right)  \nu(dy) = 0.
    \end{equation*}
    Thus,
    \begin{equation*}
        \lim_{|\xi| \to \infty} \frac{|\xi \cdot Q \xi| }{2|\xi|^2}
        \leq \lim_{|\xi| \to \infty} \frac{\re \psi(\xi)}{|\xi|^2} + \lim_{|\xi| \to \infty} \frac{1}{|\xi|^2} \lint \left(1-\cos(y \cdot \xi)\right)  \nu(dy)
        = 0
    \end{equation*}
    	which implies $Q=0$.

    \medskip\ref{app-7-ii}
    We know from \ref{app-7-i} that $Q=0$. Since $\left| \int_{|y| \geq 1} \left(1-e^{iy \cdot \xi}\right)  \nu(dy) \right| \leq 2\nu(\mbb{R}^d \setminus \ball{0}{1})$, we may assume, without loss of generality, that $\spt \nu \subseteq \cball{0}{1}$. For any $\gamma \in (0,1)$ there exists some $c_{\gamma} >0$ such that
    \begin{equation*}
        \int |y|^{\gamma} \, \nu(dy)
    	= c_{\gamma} \iint \left(1-\cos(y \cdot z)\right)\, \frac{dz}{|z|^{1+\gamma}} \, \nu(dy)
    	= c_{\gamma} \int \re \psi(z)\,\frac{dz}{|z|^{1+\gamma}}.
    \end{equation*}
    As $\spt \nu \subseteq \cball{0}{1}$, it follows easily from Taylor's formula that $|\re \psi(z)| \leq C' |z|^2$ for some absolute constant $C'>0$. On the other hand, by assumption, $|\re \psi(z)| \leq C |z|^{\beta}$ for some $\beta \in (\beta_{\infty},1)$. Consequently, we find $\lint |y|^{\gamma} \, \nu(dy)<\infty$ for all $\gamma>\beta$. This implies, in particular, that
    \begin{equation*}
        \psi_0(\xi) := \lint \left(1-e^{iy \cdot \xi}\right)  \nu(dy), \qquad \xi \in \mbb{R}^d,
    \end{equation*}
    is well-defined. Using Markov's inequality and the elementary estimate $|\sin x| \leq |x|$,  we find for all $\gamma \in (\beta,1)$
    \begin{align*}
    		|\im \psi_0(\xi)|
    		&\leq \lint_{|y \xi| < 1} |\sin(y \cdot \xi)| \, \nu(dy) + \int_{|y \cdot \xi| \geq 1} 1 \, \nu(dy)\\
    		&\leq \int_{|y \cdot \xi| <1} |y \cdot \xi| \, \nu(dy) + \int_{|y \cdot \xi| \geq 1} |y \cdot \xi|^{\gamma} \, \nu(dy)
    		\leq  |\xi|^{\gamma} \lint |y|^{\gamma} \, \nu(dy).
    \end{align*}
    Thus,
    \begin{align*}
    		2 |\xi|^{\gamma} \lint |y|^{\gamma} \, \nu(dy)
    		\geq |\im \psi_0(\xi)|
    		\geq \left| b + \lint_{|y| \leq 1} y \, \nu(dy) \right| \, |\xi| - |\im \psi(\xi)|.
    \end{align*}
    Dividing both sides by $|\xi|^{\gamma}$ and letting $|\xi| \to \infty$ proves the assertion.
\end{proof}

\defaultsection

\end{document}